\documentclass{amsart}
\usepackage{tabu}
\usepackage[margin=95pt]{geometry}
\usepackage{enumitem}
\allowdisplaybreaks

\usepackage{amssymb}
\usepackage{amsmath}
\usepackage{amscd}
\usepackage{amsbsy}
\usepackage{comment}
\usepackage{enumitem}
\usepackage[matrix,arrow]{xy}
\usepackage{hyperref}
\usepackage{xurl}
\usepackage{caption}  
\DeclareMathOperator{\Cl}{Cl}
\DeclareMathOperator{\Norm}{Norm}

\DeclareMathOperator{\ord}{ord}

\newcommand{\Q}{{\mathbb Q}}
\newcommand{\Z}{{\mathbb Z}}

\newcommand{\F}{{\mathbb F}}

\newcommand{\OO}{{\mathcal O}}

\newcommand{\gn}{{\mathfrak{n}}}

\newcommand{\fq}{\mathfrak{q}}
\newcommand{\fP}{\mathfrak{P}}

\newcommand{\fA}{\mathfrak{A}}
\newcommand{\fC}{\mathfrak{C}}
\newcommand{\fB}{\mathfrak{B}}

\def\mod#1{{\ifmmode\text{\rm\ (mod~$#1$)}
\else\discretionary{}{}{\hbox{ }}\rm(mod~$#1$)\fi}}

\begin {document}

\newtheorem{thm}{Theorem}
\newtheorem{lem}{Lemma}[section]
\newtheorem{prop}[lem]{Proposition}

\theoremstyle{definition}

\theoremstyle{remark}

\title[The Lebesgue-Nagell equation]
{Differences between perfect powers : the Lebesgue-Nagell equation}
%\author{Samir Siksek}
%\address{Mathematics Institute\\
%	University of Warwick\\
%	Coventry\\
%	CV4 7AL \\
%	United Kingdom}

%\email{s.siksek@warwick.ac.uk}
\author[Michael Bennett]{Michael A. Bennett}
\address{Department of Mathematics, University of British Columbia, Vancouver, B.C., V6T 1Z2 Canada}
\email{bennett@math.ubc.ca}

\author{Samir Siksek}
\address{Mathematics Institute, University of Warwick, Coventry CV4 7AL, United Kingdom}
\email{S.Siksek@warwick.ac.uk}
\thanks{The first-named author is supported by NSERC. The second-named author is
supported by EPSRC Grant EP/S031537/1 \lq\lq Moduli of
elliptic curves and classical Diophantine problems\rq\rq.}

\date{\today}


\begin {abstract}
We develop a variety of new techniques to treat Diophantine equations of the shape $x^2+D =y^n$, based upon bounds for linear forms in $p$-adic and complex logarithms, the modularity of Galois representations attached to Frey-Hellegouarch elliptic curves, and machinery from Diophantine approximation.
We use these to explicitly determine the set of all coprime integers $x$ and $y$,  and $n \geq 3$, with the property that $y^n > x^2$ and $x^2-y^n$ has no prime divisor exceeding $11$.
\end {abstract}
\maketitle

%------------------------------
\section{Introduction}
%------------------------------

Understanding the gaps in the sequence of positive perfect powers
$$
1, 4, 8, 9, 16, 25, 27, 32, 36, 49, 64, 81, \ldots
$$
is a problem at once classical and fundamentally difficult. Mih\u{a}ilescu's Theorem \cite{Mih} (n\'ee Catalan's Conjecture) tells us
 that $8$ and $9$ are the only consecutive integers here, but it is not, for instance, a consequence of current technology that there are at most finitely many gaps of length $k$, for any fixed integer $k > 1$ (though this was conjectured to be the case by Pillai; see e.g. \cite{Pillai}). If we simplify matters by considering instead gaps between squares and other perfect powers, then we can show that such gaps, if nonzero, grow as we progress along the sequence. Indeed, the same is even true of the greatest prime factor of the gaps. Specifically, we have the following, a special case of  Theorem 2 of Bugeaud \cite{Bug-BLMS}; here, by $P(m)$ we denote the greatest prime divisor of a nonzero integer $m$.

  \begin{thm}[Bugeaud]
  Let $n \geq 3$ be an integer. There exists an effectively computable positive constant $c=c(n)$ such that if $x$ and $y$ are coprime positive integers with $y \geq 2$, then
  $$
  P(x^2-y^n)\geq c   \log n
  $$
 and, for suitably large $x$,
  $$
  P(x^2-y^n) \geq \frac{\log \log y}{30n}.
  $$
  \end{thm}

This result is a consequence of bounds for linear forms in logarithms, complex and $p$-adic. As such, it can be made completely explicit and leads to an algorithm for solving the {\it Lebesgue-Nagell equation}
\begin{equation} \label{LebNag}
x^2+D = y^n,
\end{equation}
where we suppose that $x$ and $y$ are coprime nonzero integers, and that  either
\begin{enumerate}[label=(\roman*)]
\item $D$ is a fixed integer, or
\item all the prime divisors of $D$ belong to a fixed set of primes $S$.
\end{enumerate}
The terminology here stems from the fact that equation (\ref{LebNag}) with $D=1$
was first solved by V.\ A.\ Lebesgue \cite{Leb}, while T.\ Nagell \cite{Nag1},
\cite{Nag3} was the first researcher to study such equations in a systematic
fashion.

Regrettably, this algorithm is still, in most instances, not a practical one.
Even in the very special case $D=-2$, we are not able to completely solve
equation (\ref{LebNag}) (though there are a number of partial results available
in the literature; see e.g.\ Chen \cite{Che}). Almost all the (very ample) literature on this problem concerns cases where $D > 0$ and $y$ is odd in (\ref{LebNag}). Under these assumptions, we may solve the equation through appeal to a beautiful result of Bilu, Hanrot and Voutier \cite{BHV} on primitive divisors in binary recurrence sequences, at least for all but a few small values of $n$. Proposition 5.1 of \cite{BMS2} (sharpening work of Cohn \cite{Coh}) provides a very explicit summary of this approach -- one bounds the exponent $n$ in (\ref{LebNag}) in terms of the class numbers of a finite collection of imaginary quadratic fields, depending only upon the primes dividing $D$; see Section \ref{PrimDiv} for details. Smaller values of $n$ may be treated via techniques from elementary or algebraic number theory, or through machinery from Diophantine approximation.
By way of example, in cases (i) and (ii), equation (\ref{LebNag}), for fixed $n$, reduces to finitely many {\it Thue} or {\it Thue-Mahler} equations, respectively. These can be solved through arguments of Tzanakis and de Weger \cite{TW},  \cite{TW1}, \cite{TW2} (see also \cite{GKMS} for recent refinements).

In case either $D > 0$ and $y$ is even, or if $D < 0$, 
 the literature on equation (\ref{LebNag})  is much sparser, primarily since the machinery of primitive divisors is no longer applicable. In these cases, other than bounds for linear forms in logarithms, the only general results that we know to apply to equation (\ref{LebNag}) are derived from the modularity of Galois representations arising from associated Frey-Hellegouarch curves. These are obtained by viewing (\ref{LebNag}) as a ternary equation of signature $(n,n,2)$, i.e. as $y^n - D \cdot 1^n = x^2$. Such an approach can work to solve equation (\ref{LebNag}) in one of two ways, either by

\begin{enumerate}[label=(\alph*)]
\item producing an upper bound upon $n$ that is sharper than that coming from linear forms in logarithms, leaving a feasible set of small $n$ to treat, or
\item failing to produce such an upper bound, but, instead, providing additional arithmetic information that allows one to solve all the remaining Thue or Thue-Mahler equations below the bound coming from linear forms in logarithms.
\end{enumerate}
An example of situation (a) is the case where $D$ is divisible by only the primes in $S= \{ 5, 11 \}$ and $y$ is even. Then Theorem 1.5 of \cite{BenS} implies that equation (\ref{LebNag}) has no nontrivial solutions for all prime $n > 11$ and $y$ even;  work of Soydan and Tzanakis \cite{SoTz} treats smaller values of $n$ and the case where $y$ is odd (where the Primitive Divisor Theorem works readily). In general, we are potentially in situation (a) precisely when there fails to exist an elliptic curve $E/\mathbb{Q}$ with nontrivial rational $2$-torsion and conductor
$$
N_{S^*} = 2 \prod_{p \in S^*} p,
$$
for each subset $S^* \subseteq S$ with the property that  the product of the primes in $S^*$ is congruent to $-1$ modulo $8$. Other examples of such sets $S$ include
$$
\{ 5, 19 \}, \;  \{ 11, 13, 41 \}, \;   \{ 11, 17, 29 \}, \; \{ 11, 17, 37 \}, \; \{ 17, 19, 37 \}   \; \mbox{ and } \; \{ 19, 29 \}.
$$

For situation (b), papers of Bugeaud, Mignotte and the second author \cite{BMS2}, and of Barros \cite{Barr} deal with a number of cases of equation (\ref{LebNag}) with $D$ fixed and positive or negative, respectively.

In this  paper, we will concentrate on the first of the two difficult cases, namely when $D > 0$ and $y$ is even in (\ref{LebNag}) (so that necessarily $D \equiv -1 \mod{8}$), under the additional hypothesis that $D$ is divisible only by a few small primes. For completeness, we will also treat the easier situation where $y$ is odd, under like hypotheses on $D$. In a companion paper \cite{BeSiNew}, we will consider equation (\ref{LebNag}) in the other challenging situation where $D < 0$.
Our main result  in the paper at hand  is the complete resolution of equation (\ref{LebNag}) in case $D > 0$,  $P(D) < 13$, $\gcd (x,y)=1$  and $n \geq 3$. We prove the following.

\begin{thm} \label{Main-thm}
There are precisely $1240$ triples of positive integers $(x, y, n)$  with $n \geq 3$, $\gcd (x,y)=1$, $y^n > x^2$ and 
$$
P(x^2-y^n) < 13.
$$
They are distributed as follows.
$$
\begin{array}{|c|c|c|c|c|c|c|c|} \hline
n & \# (x,y) & n & \# (x,y)  & n & \# (x,y) & n & \# (x,y)  \\ \hline
3 & 755 & 7 & 5 & 12 & 4 & 26 & 1 \\
4 & 385 & 8 & 17 & 13 & 1  & &  \\
5 & 11 & 9 & 1 & 14 & 4  & & \\
6 & 51 & 10 & 4 &  15 &  1 & & \\ \hline
\end{array}
$$
\end{thm}
We
provide the complete list of the $1240$ solutions
at
\begin{center}
\noindent
\url{http://homepages.warwick.ac.uk/staff/S.Siksek/progs/lebnag/lebesgue_nagell_solutions.txt}
\end{center}
Proving this result amounts to solving the equation
\begin{equation} \label{eq-main}
x^2+2^{\alpha_2} 3^{\alpha_3}5^{\alpha_5}7^{\alpha_7}11^{\alpha_{11}} = y^n,
\end{equation} 
where $x, y$ and $n$ are positive integers, with $\gcd(x,y)=1$, $n \geq 3$, and the $\alpha_i$ are nonnegative integers, i.e. equation (\ref{LebNag}), where $D > 0$ is supported only on primes in $S= \{ 2, 3, 5, 7, 11 \}$. We note that earlier work along these lines typically  either treat cases where there are no $S$-units congruent to $-1 \mod{8}$, so that the analogous equations cannot have $y$ even (see e.g. the paper of Luca \cite{Lu} for $S=\{ 2, 3 \}$), or simply exclude these cases (see Pink \cite{Pink} for $S=\{2, 3, 5, 7 \}$, where solutions with $y$ even are termed {\it exceptional}). The only exceptions to this in the literature, of which we are aware, are the  aforementioned paper of Soydan and Tzanakis \cite{SoTz} where $S=\{ 5, 11 \}$ and work of Koutsianas \cite{Kou} treating $S=\{ 7 \}$ with prime exponent $n \equiv 13, 23 \mod{24}$. A comprehensive survey of the extensive literature on this equation can be found in the paper of Le and Soydan \cite{LeSo}.

To solve equation (\ref{eq-main}) completely, we are forced to introduce a variety of new techniques, many of which are applicable in rather more general settings.
These include
\begin{itemize}
\item appeal to bounds for linear forms in two $p$-adic logarithms; what is interesting here is that the resulting inequalities are surprisingly strong, leading to problems involving complex logarithms that are essentially the same level of difficulty as for the apparently easier case where $D>0$ is fixed in equation (\ref{LebNag}) (as treated in \cite{BMS2})
\item efficient sieving with Frey-Hellegouarch curves; on some level, this is likely the most important computational innovation in this paper
\item refined use of lower bounds for  linear forms in two and three complex logarithms
\item a computationally efficient approach to treat the genus one curves encountered when solving equation (\ref{eq-main}) for $n \in \{ 3, 4 \}$
\item new practical techniques for solving Thue-Mahler equations of moderate ($n \leq 13$) degree.
\end{itemize}

The outline of this paper is as follows. In Section \ref{34}, we deal with the cases of exponents $3$ and $4$ in equation (\ref{eq-main}).
In Section \ref{PrimDiv}, we apply the Primitive Divisor Theorem to handle larger exponents in (\ref{eq-main}), under the assumption that the variable $y$ is odd. 
Section \ref{TM!} begins our treatment of the complementary (significantly
harder) situation when $y$ is even, showing how equation (\ref{eq-main}) with
fixed exponent $n$ reduces to solving a number of Thue-Mahler equations. From
this, we are able to solve (\ref{eq-main}) completely for $n \leq 11$.  %changed 13 to 11
In
Section \ref{FreyH}, we show how to associate to a putative solution of
(\ref{eq-main})  a Frey-Hellegouarch elliptic curve. We then use this
connection  to develop a number of computational sieves that enable us to show
that equation (\ref{eq-main}) has no solutions with prime exponents $n$ between
$17$ and a reasonably large upper bound (which depends upon $D$, but is, in all
cases, of order exceeding $10^8$). 
This approach also deals with $n=13$, 
except for one case that is solved through reduction to a Thue--Mahler equation.
Finally, in Section \ref{sec:large}, we apply inequalities for linear forms in $p$-adic and complex logarithms to show that (\ref{eq-main}) has no solutions for exponents $n$ exceeding these upper bounds.

%---------------------------------
\section{(Very) small values of $n$} \label{34}
%----------------------------------

We begin by treating equation (\ref{eq-main}) in case $n \in \{ 3, 4 \}$. With these handled, we will thus be able to assume, without loss of generality,  that $n \geq 5$ is prime. It is worth observing that our methods of proof in this section work equally well in the analogous situation where $D$ is supported on $S= \{ 2, 3, 5, 7, 11 \}$, but $D < 0$ (a conclusion that is far from true regarding our techniques for handling larger exponents).

\subsection{Exponent $n=3$}

If we suppose that  $n=3$ in equation (\ref{eq-main}), then the problem reduces to one of determining $S$-integral points  on  
$$
3^{\#S} = 3^5=243
$$
 {\it Mordell} elliptic curves of the shape $y^2=x^3-k$, where 
$$
k = 2^{\delta_2} 3^{\delta_3}5^{\delta_5}7^{\delta_7}11^{\delta_{11}}, \; \; \mbox{ for } \; \; \delta_p \in \{ 0, 1, 2 \}.
$$
There are various ways to carry this out; if we try to do this directly using, say, the \texttt{Magma} computer algebra package \cite{magma}, we very quickly run into problems arising from the difficulty of unconditionally  certifying  Mordell-Weil bases for some of the corresponding curves. We will instead argue somewhat differently.

Given a solution to equation (\ref{eq-main}) in coprime integers $x$ and $y$, consider the Frey-Hellegouarch elliptic curve
$$
E_{x,y} \; \; : \; \; Y^2 = X^3 -3 y  X + 2 x, 
$$
with corresponding  discriminant
% \margnote{n.b. changed my mind about the minus sign!}
$$
\Delta_{E_{x,y}} \, = \,  2^{\alpha_2+6} 3^{\alpha_3+3}5^{\alpha_5}7^{\alpha_7}11^{\alpha_{11}}.
$$
This model has $c$-invariants
$$
 c_4 = 144 y \;  \mbox{ and }  \; c_6 = -1728 x.
$$
We may check via Tate's algorithm that this curve is minimal at all primes $p \geq 3$ and, while possibly not minimal at $2$, 
the fact that $x$ and $y$ are coprime implies that a corresponding  minimal model over $\mathbb{Q}$ has either 
$$
c_4 = 144 y, \; c_6 = -1728 x \; \mbox{ or } \;  c_4 = 9 y, \; c_6 = -27 x,
$$
with the latter case occurring only if $xy$ is odd.
%t is (see e.g. Papadopolous) minimal at $5, 7$ and $11$, where we find that it has either multiplicative or good reduction.
%At the prime $3$, the curve is also minimal and we have a conductor divisible by precisely $3^2, 3^3$ or $3^4$, in all cases

The isomorphism classes of elliptic curves over $\mathbb{Q}$ with good reduction outside $\{ 2, 3, 5, 7, 11 \}$ have recently been completely and rigorously determined using two independent approaches, by von Kanel and Matschke \cite{KM} (via computation of $S$-integral points on elliptic curves, based upon bounds for elliptic logarithms), and by the first author, Gherga and Rechnitzer \cite{BeGhRe} (using classical invariant theory to efficiently reduce the problem to solutions of cubic Thue-Mahler equations). One finds that there are precisely
 $592192$ isomorphism classes of  elliptic curves over $\mathbb{Q}$ with good reduction outside $\{ 2, 3, 5, 7, 11 \}$; details are available at, e.g. 
 \begin{center}
 \noindent
 \url{https://github.com/bmatschke/s-unit-equations/blob/master/elliptic-curve-tables/good-reduction-away-from-first-primes/K_deg_1/curves_K_1.1.1.1_S_2_3_5_7_11.txt}
 \end{center}
 
  For each such class, we consider the corresponding $c$-invariants;  if both $c_4 \equiv 0 \mod{144}$ and $c_6 \equiv 0 \mod{1728},$ we define
\begin{equation} \label{xy1}
y = \frac{c_4}{144} \mbox{ and } x= \frac{|c_6|}{1728},
\end{equation}
while if at least one of $c_4 \equiv 0 \mod{144}$ or $c_6 \equiv 0 \mod{1728}$ fails to hold, but  we have  $c_4 \equiv 0 \mod{9}$ and $c_6 \equiv 0 \mod{27}$, we define
\begin{equation} \label{xy2}
y = \frac{c_4}{9} \mbox{ and } x= \frac{|c_6|}{27}.
\end{equation}
For the resulting pairs $(x,y)$, we check that $y > 0$ and  $\gcd (x,y)=1$. We find $755$ such pairs, corresponding to $812$ triples $(x,y,n)$ satisfying $3 \mid n$. There are $5$ triples with $y > 10^9$, with the largest value of $y$ corresponding to the identity
$$
280213436582801^2 + 2^{16} \cdot 3^6 \cdot 5 \cdot 7^8 \cdot 11^2 = 4282124641^3.
$$

%---------------------------------
\subsection{Exponent $n=4$}
%------------------------------

In this case, we may rewrite equation (\ref{eq-main}) as
$$
(y^2-x)(y^2+x)= 2^{\alpha_2} 3^{\alpha_3}5^{\alpha_5}7^{\alpha_7}11^{\alpha_{11}} 
$$
and so either $\alpha_2=0$, in which case
\begin{equation} \label{unit1}
u_1+u_2=2y^2,
\end{equation}
where $u_i$ are coprime $\{ 3, 5, 7, 11 \}$-units, or we have
\begin{equation} \label{unit2}
u_1+u_2=y^2,
\end{equation}
where $u_i$ are coprime $\{ 2, 3, 5, 7, 11 \}$-units. In each case, since $xy \neq 0$, we may suppose that $u_1 > u_2$. To be precise, we have  
$$
u_1 u_2 =  3^{\alpha_3}5^{\alpha_5}7^{\alpha_7}11^{\alpha_{11}}, \; \; \sqrt{\frac{1}{2} (u_1+u_2)}=y \; \; \mbox{ and } \; \; 
\frac{1}{2} (u_1-u_2) = x,
$$
and 
$$
u_1 u_2 = 2^{\alpha_2-2} 3^{\alpha_3}5^{\alpha_5}7^{\alpha_7}11^{\alpha_{11}}, \; \; \sqrt{u_1+u_2}=y \; \; \mbox{ and } \; \; 
u_1-u_2 = x,
$$
in cases (\ref{unit1}) and (\ref{unit2}), respectively.

As for $n=3$, we can write down corresponding Frey-Hellegouarch curves which have good reduction outside $\{ 2, 3, 5, 7, 11 \}$ (and, additionally, in this situation, have nontrivial rational $2$-torsion). It is easier to attack this problem more directly. Both equations (\ref{unit1}) and (\ref{unit2}) take the form $a+b=c^2$, where $a$ and $b$ are $\{ 2, 3, 5, 7, 11 \}$-units with $\gcd(a,b)$ square-free. Machinery for solving such problems has been developed by de Weger \cite{Weg}, \cite{Weg2}.  Data from an implementation of this by von Kanel and Matschke \cite{KM} is available at
 \begin{center}
 \noindent
 \url{https://github.com/bmatschke/solving-classical-diophantine-equations/blob/master/sums-of-units-equations/sumsOfUnitsBeingASquare__S_2_3_5_7_11.txt}
 \end{center}
 
We find that there are $1418$ pairs $(a, b)$ such that $a + b$ is a square, $\gcd(a,b)$ is square-free, $a \geq  b$, and the only primes dividing $a$ and $b$ lie in $\{ 2, 3, 5, 7, 11 \}$. 
We further restrict our attention to those with additionally $a>b \geq 1$  and either $\gcd (a,b)=1$ (in which case we take $x=a-b$, $y = \sqrt{a+b}$), or $\gcd (a,b)=2$ (whence we choose
$x=\frac{1}{2} (a-b)$ and $y=\sqrt{\frac{1}{2} (a+b)}$). This gives $385$ pairs of coprime, positive integers $x, y$ with $y^4> x^2$ and $P(y^4-x^2) < 13$. 
These pairs actually lead to $406$ triples $(x,y,n)$ with $4 \mid n$, since $17$ of the values of $y$ are squares and four of them are cubes. 
However the four cubes have already appeared in our previous computation, so altogether we obtain $402$ new triples $(x,y,n)$ with $4 \mid n$.
Together with the $812$ triples satisfying $3 \mid n$ we have altogether $1214$ triples $(x,y,n)$.
The largest  $y$ with $n=4$ corresponds to the identity
$$
1070528159^2 \, + \, 2^{18} \cdot 3^3 \cdot 5 \cdot 7^4 \cdot 11^2 \, =\, 32719^4.
$$

For the remainder of the paper, we may therefore assume that the exponent $n$ in equation (\ref{eq-main}) is prime and $\geq 5$.

%---------------------------------
\section{Primitive divisors : equation (\ref{eq-main}) with $y$ odd} \label{PrimDiv}
%------------------------------

In this section, we treat \eqref{eq-main} under the assumption that $y$ is odd, using the celebrated Primitive Divisor Theorem of Bilu, Hanrot and Voutier \cite{BHV},
and prove the following proposition.
\begin{prop}\label{prop:yodd}
The only solutions to (\ref{eq-main}) with $n \geq 5$ prime, $\gcd(x,y)=1$ and $y$ odd correspond to the identities
\begin{gather*}
1^2+ 2 \cdot 11^2 =3^5,
\quad
241^2+2^3 \cdot 11^2  =9^5,
\quad
401^2+ 2 \cdot 5^3=11^5,\\
4201^2+ 2 \cdot 3 \cdot 5^3 \cdot 11^4=31^5 \; \; \mbox{ and } \; \; 
4443^2+2^2 \cdot 7 \cdot 11^6 = 37^5.
\end{gather*}
\end{prop}
If we consider solutions with $\gcd(x,y)=1$, $y$ odd, and $n$  divisible
by a prime $\ge 5$, then we must count one more solution corresponding to
$241^2+2^3 \cdot 11^2  =3^{10}$.  Thus our total number of 
solutions to \eqref{eq-main} for cases considered so far is $1214+6=1220$. 

\subsection{Lucas Sequences and the Primitive Divisor Theorem}
It is convenient to first introduce Lucas sequences as defined in \cite{BHV}.
A pair $(\gamma,\delta)$ of algebraic integers is called a
\textbf{Lucas pair} 
if $\gamma+\delta$ and $\gamma\delta$ are
non-zero coprime rational integers, and $\gamma/\delta$
is not a root of unity. 
%In particular, associated
%to the Lucas pair $(\alpha,\beta)$ is a
%\textbf{characteristic polynomial}
%\[
%X^2-(\alpha+\beta) X + \alpha \beta \; \in \; \Z[X].
%\]
%This polynomial has discriminant $D=(\alpha-\beta)^2 \in \Z \setminus \{0\}$.
Given a Lucas pair $(\gamma,\delta)$ we define
the corresponding \textbf{Lucas sequence} by
\[
L_m=\frac{\gamma^m-\delta^m}{\gamma-\delta}, \qquad m=0,1,2,\dotsc.
\]
A prime $\ell$ is said to be a
\textbf{primitive divisor} of the $m$-th term 
if $\ell$ divides $L_m$ but $\ell$ does not divide $(\gamma-\delta)^2 \cdot u_1 u_2 \dotsc u_{m-1}$.
\begin{thm}[Bilu, Hanrot and Voutier \cite{BHV}]\label{thm:BHV}
Let $(\gamma,\delta)$ be a Lucas pair and write $\{L_m\}$
for the corresponding Lucas sequence.
If $m \ge 30$, then $L_m$ has a primitive divisor.
Moreover, if $m \ge 11$ is prime, then $L_m$ has a primitive
divisor.
\end{thm}

Let $\ell$ be a prime. 
We define the \textbf{rank of apparition of $\ell$}
in the Lucas sequence $\{L_m\}$ to be the smallest
positive integer $m$ such that $\ell \mid L_m$. We denote the
rank of apparition of $\ell$ by $m_\ell$.
The following theorem of Carmichael \cite{Car} will be useful to us; for a concise proof
see \cite[Theorem 8]{BGPS}.
\begin{thm}[Carmichael  \cite{Car}]\label{thm:Carmichael}
Let $(\gamma,\delta)$ be a Lucas pair, and $\{L_m\}$
the corresponding Lucas sequence. Let $\ell$ be a prime.
\begin{enumerate}[label=(\roman*)]
\item If $\ell \mid \gamma \delta$ then $\ell \nmid L_m$
for all positive integers $m$.
\item Suppose $\ell \nmid \gamma \delta$.
Write $D=(\gamma-\delta)^2 \in \Z$.
\begin{enumerate}[label=(\alph*)]
\item If $\ell \ne 2$ and $\ell \mid D$, then
$m_\ell=\ell$.
\item If $\ell \ne 2$ and $\left(\frac{D}{\ell}\right)=1$, then $m_\ell \mid (\ell-1)$.
\item If $\ell \ne 2$ and $\left(\frac{D}{\ell} \right)=-1$, then $m_\ell \mid (\ell+1)$.
\item If $\ell=2$, then $m_\ell=2$ or $3$.
\end{enumerate}
\item If $\ell \nmid \gamma \delta$ then
\[
\ell \mid L_m \iff m_\ell \mid m.
\]
\end{enumerate}
\end{thm}

\subsection{Equation \eqref{eq-main} with $y$ odd}
For the remainder of this section, $(x,y,n,\alpha_2,\dotsc,\alpha_{11})$ will denote a solution to the equation 
\begin{equation}\label{eqn:yodd}
x^2+2^{\alpha_2} 3^{\alpha_3} 5^{\alpha_5} 7^{\alpha_7} 11^{\alpha_{11}} \; = \; y^n \; \; \mbox{ $x>0$, $y$ odd, } \; \gcd(x,y)=1 \; \mbox{and} \;  n \ge 5 \; \mbox{prime}.
\end{equation}
We shall write
\begin{equation} \label{fishhead}
2^{\alpha_2} 3^{\alpha_3} 5^{\alpha_5} 7^{\alpha_7} 11^{\alpha_{11}} \, = \, c^2 d, \; \;  \text{where $d$ is squarefree}.
\end{equation}
\begin{lem}\label{lem:preBHV}
There exist integers $u$ and $v$ such that
\[
x+c\sqrt{-d} \, = \, (u+v \sqrt{-d})^n, \; \mbox{ where } \;  y=u^2+d v^2, \; \;  u \mid x, \; \;  v \mid c\; \mbox{ and } \;  \gcd(u,dv)=1.
\]
If we define 
\[
\gamma=u+v\sqrt{-d} \; \; \mbox{ and } \; \;  \delta=u-v\sqrt{-d},
\]
then $(\gamma,\delta)$ is a Lucas pair. Let $\{L_m\}$ be the corresponding Lucas sequence. Then 
\begin{equation}\label{eqn:BHVTM}
L_n \, = \, \frac{\gamma^n-\delta^n}{\gamma - \delta} \, = \, \frac{c}{v}.
\end{equation}
\end{lem}
\begin{proof}
Write $M=\Q(\sqrt{-d})$.
From \eqref{eqn:yodd}, we have
\[
(x+c \sqrt{-d}) (x-c \sqrt{-d}) = y^n
\]
where the two factors on the left generate coprime ideals of $\OO_M$. Thus $(x+c \sqrt{-d}) \OO_M = \fA^n$, for some ideal of $\fA$ of $\OO_M$.
There are $32$ possible values of $d$ and
we checked, via \texttt{Magma}, that the corresponding quadratic fields $M=\Q(\sqrt{-d})$  have, in every case,  class numbers $h$ satisfying
$$
h \in \{ 1, 2, 4, 8, 12, 32 \}.
$$
In particular $n$ is coprime to $h$, and therefore $\fA$ is principal. We deduce that
$x+c\sqrt{-d}=\epsilon \cdot \gamma^n$  where $\epsilon \in \OO_M^*$ and $\gamma \in \OO_M$.
The order of the unit group $\OO_M^*$ is either $4$ (if $d=1$),  $6$ (if $d=3$) or $2$ (in all other cases).
Thus the unit group is $n$-divisible and we may absorb $\epsilon$ into the $\gamma^n$ factor to obtain
$x+c\sqrt{-d} = \gamma^n$ for some $\gamma \in \OO_M$.
We write $\delta$ for the conjugate of $\gamma$. Note that $\gamma+\delta$
is a divisor of $2x = \gamma^n +\delta^n$ and that $\gamma \delta=y$. It follows that
$\gamma+\delta$ and $\gamma \delta$ are non-zero coprime rational integers.
We claim that $\gamma/\delta$ is not a unit. If we suppose otherwise, then $(x+c\sqrt{-d})/(x-c\sqrt{-d})=(\gamma/\delta)^n$
is a unit. By coprimality of the numerator and denominator, we obtain that $x+c\sqrt{-d}$
is a unit and therefore $y=1$, a contradiction. Thus $\gamma/\delta$ is not a unit
and $(\gamma,\delta)$ is a Lucas pair. Write $\{L_m\}$ for the corresponding Lucas
sequence.

Since $\gamma \in \OO_M$, 
we have $\gamma=u+v \sqrt{-d}$ with $u$ and $v$ are integers,
or $\gamma=(u+v\sqrt{-d})/2$ where both $u$ and $v$ are odd integers.
Suppose first that we are in the latter case (whence we note that necessarily $d \equiv 3 \pmod{4}$).
Observe that $\gamma^n - \delta^n = 2 c \sqrt{-d}$, so that $v \sqrt{-d}=\gamma - \delta$ divides $2c \sqrt{-d}$. As $v$ is odd,
we deduce that $v \mid c$ and that
\[
L_n=\frac{\gamma^n-\delta^n}{\gamma-\delta} = 2 \cdot \frac{c}{v}.
\]
In particular, $L_n$ is even. 
We note that $\gamma+\delta=u$ and $\gamma \delta=(u^2+dv^2)/4=y$.
Thus the sequence $\{L_m\}$ satisfies the recurrence
\[
L_0=0, \quad L_1=1, \qquad L_{m+2}=uL_{m+1} - y L_{m}. 
\]
Using the fact that $u$ and $v$ are odd, one checks by induction that
\[
L_m \equiv 0 \pmod{2} \; \iff \; 3 \mid m.
\]
Thus, in particular, $3 \mid n$, contradicting the assumption that $n \ge 5$ is prime.
It follows that $\gamma=u+v \sqrt{-d}$ where $u$ and $v$ are integers.
Now observe that
\[
2u = (\gamma+ \delta) \mid (\gamma^n + \delta^n) = 2x \; \; \mbox{ and } \; \;  2v \sqrt{-d} = (\gamma-\delta) \mid (\gamma^n - \delta^n)=2c \sqrt{-d}.
\]
Thus $u \mid x$ and $v \mid c$. Since $y=u^2+dv^2$, we conclude that $\gcd(u,dv)=1$.
The lemma follows.
\end{proof}

\begin{lem}
Let $(x,y,n,\alpha_2,\dotsc,\alpha_{11})$ be a solution to \eqref{eqn:yodd}.
Then $n=5$ and 
%$L_n= \pm 5^r \cdot 11^s$ for integers $r \ge 0$ and $s \ge 1$.
\begin{equation} \label{eq5}
5 u^4- 10 d u^2v^2  + d^2 v^4 \,  =\,   \pm 5^r \cdot 11^s,
\end{equation}
for some $r \in \{0,1\}$ and $s \ge 0$. 
%Moreover, $r=1$ precisely when $5 \mid dv$.
\end{lem}
\begin{proof}
%By Lemma~\ref{lem:preBHV}, the only possible prime divisors of $L_n$
%are $2$, $3$, $5$, $7$ and $11$. 
We continue with the notation of Lemma~\ref{lem:preBHV}.
By \eqref{eqn:BHVTM}, we have $L_n=c/v$;
this is coprime to $\gamma \delta=y$.
If $\ell$ is any prime divisor of  $L_n$ 
then $m_\ell=n$, 
by part (iii) of Theorem~\ref{thm:Carmichael}
and the primality of $n$.

Suppose first that $n \ge 11$. 
By Theorem~\ref{thm:BHV}, $L_n$ must have a primitive divisor,
$q$ say.
By definition, this does not divide $D=(\gamma-\delta)^2=-4v^2 d$.
Thus by part (ii) of Theorem~\ref{thm:Carmichael}, $n \mid (q-1)$ if $(D/q)=1$ 
and $n \mid (q+1)$ if $(D/q)=-1$.
As the possible values of $q \mid (c/v)$ are $2$, $3$, $5$, $7$ and $11$,
we obtain a contradiction. Thus $n=5$ or $n=7$.
%and $n \ge 5$ is prime, we deduce that $n=5$, that the primitive divisor
%must be $q=11$. Moreover $11 \nmid dv$ and $(-d/11)=1$.

Next we deal with the case $n=7$. If we suppose that $L_7$
has a primitive divisor $q$ then the above argument shows
that $7 \mid (q-1)$ or $7 \mid (q+1)$ which is impossible
as $q \in \{2,3,5,7,11\}$. Thus $L_7$ has no primitive
divisor and our Lucas pair $(\gamma,\delta)=(u+v\sqrt{-d},u-v\sqrt{-d})$
is $7$-defective in the terminology of \cite{BHV}. In particular,
by the classification of defective Lucas pairs (Theorem C of \cite{BHV}) we have
$u+v\sqrt{-d}=\pm (1+\sqrt{-7})/2$ or $\pm (1+\sqrt{-19})/2$.
Both are impossible as $u$ and $v$ are integers. Hence there are no
solutions to \eqref{eqn:yodd} with $n=7$. 

Finally we deal with $n=5$.
Since $5 \nmid \ell \cdot (\ell-1)(\ell+1)$ for any of $\ell \in \{ 2, 3, 7 \}$,
we see that $L_n=c/v= \pm 5^r \cdot 11^s$. Moreover, if $r \ge 1$
then $m_5=5$ and so $5 \mid dv$ by Theorem~\ref{thm:Carmichael}.

Substituting $\gamma=u+v\sqrt{-d}$ and $\delta=u-v\sqrt{-d}$
in \eqref{eqn:BHVTM} gives \eqref{eq5}. We note that if $r \ge 2$
then $5 \mid u$, contradicting the coprimality of $u$ and $dv$.
Therefore $r \in \{0,1\}$.
\end{proof}
It remains to solve the quartic Thue--Mahler equations \eqref{eq5} for our $32$ possible
values of $d$.
Appealing to  the Thue-Mahler equation solver, implemented in  \texttt{Magma} and associated to the paper \cite{GKMS},
we obtain the following solutions:
\[
(d,u,v) \, = \, (2,-1,1), \; \;  (2,1,-2), \; \;  (7,3,-2), \; \;  (10,1,1) \; \mbox{ and } \;  (30,1,1).
\]
\begin{comment}
$$
n=5, \; \; d=2, \; \; (u,v)=(1,1), (1,2), 
$$
$$
n=5, \; \;  d=7, \; \; (u,v)=(3,2),  
$$
$$
n=5, \; \;  d=10, \; \; (u,v)=(1,1)
$$
and
$$
n=5, \; \;  d=30, \; \; (u,v)=(1,1).
$$
\end{comment}
These  lead, respectively,  to solutions of equation (\ref{eq-main}) with
$$
(x,y,n) \, =\, (1,3,5),\; (241,9,5),\; (4443,37,5),\; (401,11,5) \; \mbox{ and } \;  (4201,31,5),
$$
completing the proof of Proposition~\ref{prop:yodd}.

\section{Reduction to Thue-Mahler equations: the case of even $y$} \label{TM!}
%------------------------------------------------------------------------

From the results of the preceding sections, we are left to treat \eqref{eq-main} with $y$ even and $n \ge 5$ prime. 
It therefore remains to consider the equation
\begin{equation}\label{eqn:yeven}
x^2+3^{\alpha_3} 5^{\alpha_5} 7^{\alpha_7} 11^{\alpha_{11}} \; = \; y^n \; \; \mbox{ with $y$ even, } \; \gcd(x,y)=1 \; \mbox{ and } \;  n \ge 5 \; \mbox{ prime}.
\end{equation}
Let us define
\begin{equation} \label{n-upper-bound}
N(d) =
\begin{cases}
 6 \times 10^8 \; \mbox{ if } \; d=7, \\
 4 \times 10^8 \; \mbox{ if } \; d=15, \\
 5 \times 10^8 \; \mbox{ if } \; d=55, \\
 1.2 \times 10^9 \; \mbox{ if } \; d=231. \\
 \end{cases}
 \end{equation}
The purpose of this section and the next is to prove the following proposition.
\begin{prop}\label{prop:yeven}
The only solutions to \eqref{eqn:yeven} with $n<N(d)$ correspond to the identities
\begin{gather*}
31^2+ 3^2 \cdot 7 =4^5,
\quad
5^2+7  =2^5,
\quad
181^2+ 7=8^5,
\quad
17^2+ 3\cdot 5 \cdot 7^2 = 4^5,\\
23^2+ 3^2 \cdot 5 \cdot 11=4^5,
\quad
130679^2+3 \cdot 7^3 \cdot 11^7 = 130^5,
\quad
47^2 + 3^4\cdot 5^2\cdot 7=4^7,
\quad
	11^2+7 = 2^{7},
	\\
7^2 + 3^3 \cdot 5 \cdot 11^2=4^7,
\quad
117^2+ 5 \cdot 7^2 \cdot 11 =4^7,
\quad
103^2+3 \cdot 5^2 \cdot 7 \cdot 11 = 4^7,
\\
\mbox{ and } \; \;   8143^2+3^3 \cdot 5 \cdot 7^2 \cdot 11^2 = 4^{13}.
\end{gather*}
\end{prop}
This gives $12$ new solutions to \eqref{eq-main}
with $n \in \{ 5, 7, 13 \}$ and, additionally,  
$8$ further solutions with exponents $10$, $14$ and $26$.
Thus the total number of solutions we have found
so far for \eqref{eq-main} is $1220+12+8=1240$.
We shall show in Section~\ref{sec:large}
that there are no further solutions,
and that therefore \eqref{eq-main} has precisely
$1240$ solutions as claimed in Theorem~\ref{Main-thm}.
	\begin{comment}
    [ 5, 31, 4, 0, 2, 0, 1, 0 ],
    [ 5, 5, 2, 0, 0, 0, 1, 0 ],
    [ 5, 181, 8, 0, 0, 0, 1, 0 ],
    [ 5, 17, 4, 0, 1, 1, 2, 0 ],
    [ 5, 23, 4, 0, 2, 1, 0, 1 ],
    [ 5, 130679, 130, 0, 1, 0, 3, 7 ]
]
,
[
    [ 7, 47, 4, 0, 4, 2, 1, 0 ],
    [ 7, 11, 2, 0, 0, 0, 1, 0 ],
    [ 7, 7, 4, 0, 3, 1, 0, 2 ],
    [ 7, 117, 4, 0, 0, 1, 2, 1 ],
    [ 7, 103, 4, 0, 1, 2, 1, 1 ]
		[ 13, 8143, 4, 0, 3, 1, 2, 2 ]
	\end{comment}

We assume without loss of generality that $x \equiv 1 \mod{4}$.
As before we shall write
\begin{equation}\label{eqn:sqfree}
3^{\alpha_3} 5^{\alpha_5} 7^{\alpha_7} 11^{\alpha_{11}} \; = \;  c^2 d, \;  \text{where $d$ is squarefree and }   c \; = \; 3^{ \beta_3} 5^{\beta_5} 7^{\beta_7} 11^{ \beta_{11}}.
\end{equation}
Since $y$ is even, it follows from \eqref{eqn:yeven} that $d \equiv - 1 \mod{8}$, whence necessarily
\begin{equation}\label{eqn:dvals}
d \in \{ 7, 15, 55, 231 \}.
\end{equation}
Let $M=M_d=\Q(\sqrt{-d})$.
We note the structure of the class group of $M$:
\[
\Cl(M)  \cong \begin{cases}
1 & d=7\\
C_2 & d=15\\
C_4 & d=55\\
C_2 \times C_6 & d=231.
\end{cases}
\]

%[ 7, 15, 55, 231 ]
%Abelian Group of order 1
%7 1 1/4*(-th + 1)
%Abelian Group isomorphic to Z/2
%Defined on 1 generator
%Relations:
%    2*Cl.1 = 0
%15 2 1/8*(-th + 1)
%Abelian Group isomorphic to Z/4
%Defined on 1 generator
%Relations:
%    4*Cl.1 = 0
%55 4 1/32*(-th - 3)
%Abelian Group isomorphic to Z/2 + Z/6
%Defined on 2 generators
%Relations:
%    2*Cl.1 = 0
%    6*Cl.2 = 0
%231 6 1/128*(-th + 5)

\begin{lem}\label{lem:facteven}
Let $c^\prime = \pm c$ with the sign chosen so that $c^\prime \equiv 1 \mod{4}$.
Let
\[
h=\begin{cases}
1 & d=7\\
2 & d=2\\
4 & d=55\\
6 & d=231
\end{cases}
\; \mbox{ and } \; 
\eta=r+s \sqrt{-d}, \; \mbox{ where } \; 
(r,s) \; = \; 
\begin{cases}
(1/4,-1/4) & d=7\\
(1/8,-1/8) & d=15\\
(3/32,1/32) & d=55\\
(5/128,-1/128) & d=231.
\end{cases}
\]
Let $0 \le \kappa_n \le h-1$ be the unique integer satisfying
$\kappa_n \cdot n \equiv -2 \mod{h}$. 
Then there is some non-zero $\mu \in \OO_M$ such that
\begin{equation}\label{eqn:principal1}
\frac{x +c^\prime  \sqrt{-d}}{2}  \; = \;
\eta^{(2+\kappa_n \cdot n)/h} \cdot \mu^n.
\end{equation}
Moreover, $\eta$ is supported only on prime ideals dividing $2$
and $\mu$ is supported only on prime ideals dividing $y$.
\end{lem}
\begin{proof}
As $d \equiv -1 \mod{8}$, the prime $2$
splits in $\OO_M$ as $2\OO_M=\fP \cdot \overline{\fP}$,
where
\begin{equation}\label{eqn:fP}
\fP \; = \; 2  \OO_M \; + \; \left(\frac{1+\sqrt{-d}}{2} \right) \cdot \OO_M.
\end{equation}
We may rewrite \eqref{eqn:yeven} as
\begin{equation}\label{eqn:rewrite}
\left( \frac{x+c^\prime \sqrt{-d}}{2} \right)
\left( \frac{x-c^\prime \sqrt{-d}}{2} \right)
\; = \; \frac{y^n}{4}.
\end{equation}
Note that the two factors on the left hand-side of this last equation are coprime elements of $\OO_M$. 
Since $x \equiv c^\prime \equiv 1\mod{4}$, 
we see that $\fP$ divides the first factor on the left-hand-side.
We thus deduce that
\begin{equation}\label{eqn:fpfa}
\left( \frac{x +c^\prime \sqrt{-d}}{2} \right) \cdot \OO_M \; = \; \fP^{-2} \cdot \fA^n,
\end{equation}
where $\fA$ is an integral ideal divisible by $\fP$, with $\fA \cdot \overline{\fA} = y \OO_M$. 
The order of the class $[\fP]$ in $\Cl(M)$ is $h$. Thus $\fP^{-h}$ is principal,
and $\eta$ has been chosen so that $\fP^{-h}=\eta \OO_M$.
Let $\fB=\fP^{\kappa_n} \cdot \fA$. Then we may rewrite \eqref{eqn:fpfa}
as 
\[
\left( \frac{x +c^\prime \sqrt{-d}}{2} \right) \cdot \OO_M \; = \; \fP^{-(2+\kappa_n \cdot n)} \cdot \fB^n \; = \; \eta^{(2+\kappa_n \cdot n)/h } \cdot \fB^n.
\]
Since $n$ is a prime that does not divide the order of $\Cl(M)$, the ideal $\fB$ must be principal. Let $\mu$
be a generator for $\fB$. Then
\[
 \frac{x +c^\prime \sqrt{-d}}{2}   \; =  \; \pm \eta^{(2+\kappa_n \cdot n)/h } \cdot \mu^n
\]
and  \eqref{eqn:principal1} follows on absorbing the $\pm$ sign into $\mu$.
It is clear that $\eta$ is supported on $\fP$ only. 
Moreover $\fB$ is an integral ideal with norm $2^{\kappa_n} y$. It follows, since $y$ is even, that $\mu$ is supported
only on prime ideals dividing $y$.
\end{proof}

\begin{lem}\label{lem:TM}
The only solutions to equation (\ref{eqn:yeven})
with  $n \in \{ 5, 7, 11 \}$  are those corresponding to the identities in 
	Proposition~\ref{prop:yeven}.
\end{lem}
\begin{proof}
We drop our requirement that $x>0$ and replace it with the assumption $x \equiv 1 \mod{4}$, so that
we can apply Lemma~\ref{lem:facteven}. 
For each exponent $n$, there are four cases to consider depending on the value
of $d \in \{7,15,55,231\}$ in \eqref{eqn:sqfree}.
For each pair $(n,d)$, Lemma~\ref{lem:facteven}
asserts that $(x,c^\prime)$ satisfies \eqref{eqn:principal1}
with $\mu \in \OO_M$. We write
\[
\mu=r+s (1+\sqrt{-d})/2, 
\]
with $r$ and $s$ rational integers. 
We will show that $\gcd(r,s)=1$.
	If $2 \mid r$ and $2 \mid s$ then $\overline{\fP} \mid \mu$
	which contradicts the coprimality of the two factors in the
	left hand-side of \eqref{eqn:rewrite}. If $\ell$ is an odd prime with $\ell \mid r$ and $\ell \mid s$,
	then again we contradict the coprimality of those two factors.
	Hence $\gcd(r,s)=1$.

From \eqref{eqn:principal1},  we have
\[
	c^\prime = \frac{1}{\sqrt{-d}} \left( \eta^m \cdot (r+s(1+\sqrt{-d})/2)^n \, - \, \overline{\eta}^m \cdot (r+s(1-\sqrt{-d})/2)^n \right)
\]
	where $m=(2+\kappa_n \cdot n)/h$. The expression on the right has the form $2^{-hm} F(r,s)$
	where $F \in \Z[X,Y]$ is homogeneous of degree $n$. 
We therefore, in each case,  obtain a Thue-Mahler equation of the form 
\[
	F(r,s) \; =\; 2^{hm} \cdot c^\prime\; =\; 
	\pm 2^{hm} \cdot 3^{\beta_3} 5^{\beta_5} 7^{\beta_7} 11^{\beta_{11}}.
\]
We solved these Thue-Mahler equations using 
the Thue-Mahler solver associated
with the paper
\cite{GKMS}. This computation took around one day and resulted in the solutions
in Proposition~\ref{prop:yeven} for $n \in \{ 5, 7 \}$; there were no solutions
for $n=11$.
\end{proof}

%-----------------------------------------------------------------------
\section{Frey-Hellegouarch curves and related objects} \label{FreyH}
%-----------------------------------------------------------------------

We continue to treat \eqref{eq-main} with $y$ even, i.e. equation (\ref{eqn:yeven}), where we maintain the assumption that $x \equiv 1 \mod{4}$.
Although the results of the previous section allow us to assume more, for now we merely
impose the following constraint on the exponent:
	$n \ge 7$ is prime.
Following the first author and Skinner \cite{BenS}, we associate to a solution $(x,y,n)$ the  Frey-Hellegouarch elliptic curve $F=F(x,y,n)$ defined via
\begin{equation}\label{eqn:evenFrey}
F \; \; : \; \; Y^2+XY = X^3+\left( \frac{x-1}{4} \right) X^2 + \frac{y^n}{64} X.
\end{equation}
The model here is minimal, semistable, and we note the following invariants,
\[
c_4=x^2 - \frac{3}{4} y^n, \qquad c_6=-x^3+\frac{9}{8} x y^n
\]
and
\[
\Delta_F \; = \; \frac{y^{2n}}{2^{12}} ( x^2-y^n) \; = \; - 2^{-12} \cdot 3^{\alpha_3} 5^{\alpha_5} 7^{\alpha_7} 11^{\alpha_{11}} \cdot y^{2n}.
\]
We invoke work of the first author and Skinner  \cite{BenS}, building on the modularity of elliptic curves over $\Q$
following Wiles and others \cite{Wiles}, \cite{BreuilConradDiamondTaylor01}, Ribet's level lowering
theorem \cite{Ribet-1990}, and the isogeny theorem of Mazur \cite{Mazur-1978}.
Write $N$ for the conductor of $E$ and let
\[
N^\prime = \frac{N}{\prod_{\stackrel{\ell \mid\mid N}{n \mid \ord_\ell(\Delta_F)}} \ell}.
\]
The results of  \cite{BenS}  assert
the existence of a weight $2$ newform $f$ of level $N^\prime$
such that 
\begin{equation}\label{eqn:compare-first}
\overline{\rho}_{F,n} \sim \overline{\rho}_{f,\frak{n}},
\end{equation}
with $\gn \mid n$ a prime ideal in the ring of integers $\OO_K$ of the Hecke eigenfield $K$ of $f$.
\begin{lem}\label{lem:condequiv}
We have $N^\prime=2 R$ where $R \mid 3 \cdot 5 \cdot 7 \cdot 11$. Moreover, for $\ell \in \{3,5,7,11\}$,
we have
\begin{equation}\label{eqn:condequiv}
\ell \nmid N^\prime \; \iff \; \alpha_\ell \equiv 0 \pmod{n} \; \iff 2 \ord_\ell(c)+\ord_\ell(d) \equiv 0 \pmod{n}
\end{equation}
where $c$ and $d$ are given in \eqref{eqn:sqfree}.
\end{lem}
\begin{proof}
Since $E$ is semistable, $N$ is squarefree, and therefore $N^\prime$ is squarefree.
Note that $\ord_2(\Delta_F)= 2n \ord_2(y)-12$. Thus $2 \mid\mid N$ and $n \nmid \ord_2(\Delta_F)$,
whereby $2 \mid\mid N^\prime$.

Next let $\ell \ge 13$. Then $\ord_\ell(\Delta_F)=2n \ord_\ell(y)$ and hence  $\ell \nmid N^\prime$.
It follows that $N^\prime=2 R$ with $R \mid 3 \cdot 5 \cdot 7 \cdot 11$.

To prove the second part of the lemma, note that, for $\ell \in \{3,5,7,11\}$,
\[
\ord_\ell(\Delta)=\alpha_\ell+2n \ord_\ell(y)= 2\ord_\ell(c)+\ord_\ell(d)+2n \ord_\ell(y).
\]
If $\alpha_\ell=0$ and $\ord_\ell(y)=0$, then $\ell \nmid N$ and so $\ell \nmid N^\prime$,
and therefore \eqref{eqn:condequiv} holds. Suppose $\alpha_\ell>0$ or $\ord_\ell(y) >0$.
Then $\ell \mid \mid N$. By the formula for $N^\prime$, we have
$\ell \nmid N^\prime$ if and only if $n \mid \ord_\ell(\Delta)$ which is
equivalent to $n \mid \alpha_\ell$. This completes the proof.
\end{proof}

Let $f$ be the weight $2$ newform of level $N^\prime$ satisfying \eqref{eqn:compare-first}.
Write
\begin{equation}\label{eqn:qexp-first}
f=\mathfrak{q}+\sum_{m=2}^\infty c_m \mathfrak{q}^m
\end{equation}
for the usual $q$-expansion of $f$. Then  $K=\Q(c_1,c_2,\ldots)$,
and the coefficients $c_i$ belong to  $\OO_K$.
\begin{lem}\label{lem:expbound}
Let $\ell \nmid N^\prime$ be a prime and 
write
\[
	\mathcal{C}_{f,\ell}^{\prime}=
	\begin{cases}
		(\ell+1)^2 - c_\ell^2 & \text{if $K=\Q$}\\
		\ell \cdot ((\ell+1)^2 - c_\ell^2)) & \text{if $K \ne \Q$}.
	\end{cases}
\]
Let $d$ be as in \eqref{eqn:sqfree}, and set
\[
	T_\ell(f)=
\begin{cases}
	\{a \in \Z \cap [-2\sqrt{\ell},2\sqrt{\ell}] \; : \; \ell+1-a \equiv 0 \mod{4}\} & \text{if $(-d/\ell)=1$} \\
	\{a \in \Z \cap [-2\sqrt{\ell},2\sqrt{\ell}] \; : \; \ell+1-a \equiv 0 \mod{2}\} & \text{if $(-d/\ell)=-1$}\\
	\hskip25ex \emptyset & \text{if $\ell \mid d$}.
\end{cases}
\]
Let
\[
\mathcal{C}_{f,\ell} \; = \;
\mathcal{C}_{f,\ell}^\prime
	\cdot \prod_{a \in T_\ell(f)} (a-c_\ell).
\]
	If $\overline{\rho}_{F,n} \sim \overline{\rho}_{f,\gn}$, then $\gn \mid \mathcal{C}_{f,\ell}$.
\end{lem}
\begin{proof}
Suppose $\ell \nmid N^\prime$ and write $N$ for the conductor of $F$. 
Suppose $\overline{\rho}_{F,n} \sim \overline{\rho}_{f,\gn}$. 
	A standard consequence \cite[Propositions 5.1, 5.2]{Siksek} of this is that
\[
	\begin{cases}
		c_\ell \equiv a_\ell(F) \mod{\gn} & \text{if $\ell \ne n$ and $\ell \nmid N$}\\
		c_\ell \equiv \pm (\ell+1) \mod{\gn} & \text{if $\ell \ne n$ and $\ell \mid N$}.
	\end{cases}
\]
Here the restriction $\ell \ne n$ is unnecessary if $K=\Q$.
It follows if $\ell \mid N$ that $\gn \mid \mathcal{C}_{f,\ell}^\prime$.
We observe that the discriminant of $F$ can written as
\[
	\Delta \; =\; (-d) \cdot (c y^n/2^6)^2. 
\]
If $\ell \mid d$ then $\ell \mid N$ and so we take $\mathcal{C}_{f,\ell}=\mathcal{C}_{f,\ell}^\prime$.

Suppose $\ell \nmid N$ and so $\ell \nmid d$. Thus  $c_\ell \equiv a_\ell(F) \mod{\gn}$.
To complete the proof it is sufficient to show that $a_\ell(F) \in T_\ell(f)$.
The model for $F$ given in \eqref{eqn:evenFrey} is isomorphic to 
\begin{equation}\label{eqn:simplifiedModel}
F \; \; : \; \; Y^2 = X^3 +xX^2+\frac{y^n}{4} X   ,
\end{equation}
and so has a point of order $2$. Thus $\ell+1-a_\ell(F)=\# F(\F_\ell) \equiv 0 \mod{2}$.
Moreover, if $(-d/\ell)=1$ then the discriminant is a square modulo $\ell$,
so $F/\F_\ell$ has full $2$-torsion, whence 
$\#F(\F_\ell) \equiv 0 \mod{4}$. It follows that $a_\ell(F) \in T_\ell(f)$.
\end{proof}

There are a total of $76$ conjugacy classes of newforms $f$ at 
the levels $N^\prime=2R$ with $R \mid 3 \cdot 5 \cdot 7 \cdot 11$,
of which $59$ are rational (and so correspond to elliptic curves).
Since there  are four possible values of $d \in \{7,15,55,231\}$,
this gives $4 \times 76=304$ pairs $(f,d)$ to consider.
We apply Lemma~\ref{lem:expbound} to each pair $(f,d)$,
letting
\[
	\mathcal{C}_{f,d}=\sum \mathcal{C}_{f,\ell} \cdot \OO_K,
\]
where the sum is over all primes $3 \le \ell < 500$ not dividing $N^\prime$.
It follows from Lemma~\ref{lem:expbound} that $\gn \mid \mathcal{C}_{f,d}$.
We let
\[
	C_{f,d}=\Norm_{K/\Q}(\mathcal{C}_{f,d}).
\]
Since $\gn \mid n$, we have that $n \mid C_{f,d}$.
Of the $304$ pairs $(f,d)$, the integer ${C}_{f,d}$ is identically 
zero for $114$ pairs, and non-zero for the remaining $190$ pairs.
For the $190$ pairs $(f,d)$ where ${C}_{f,d} \ne 0$,
we find that the largest possible prime divisor
of any of these $C_{f,d}$ is $11$. By the results of the 
previous section we know all the solutions to
\eqref{eqn:yeven} with $n \in \{ 7, 11 \}$ and hence can
therefore eliminate these $190$ pairs from 
further consideration. We focus on the $114$ 
remaining pairs $(f,d)$.
Here, each $f$ satisfies $K=\Q$
and so corresponds
to an elliptic curve $E/\Q$ whose conductor
is equal to the level $N^\prime$
of $f$. Moreover,
each of these elliptic curve $E$ has 
non-trivial rational $2$-torsion.
This is unsurprising in view of 
the remarks following \cite[Proposition 9.1]{Siksek}.
We observe that $\overline{\rho}_{f,n} \sim \overline{\rho}_{E,n}$.
Thus we have $114$ pairs $(E,d)$ to consider,
and if $(x,y,n)$ is a solution to \eqref{eqn:yeven}
with $n \ge 13$ prime then there is some
pair $(E,d)$ (among the $114$) where $d$ satisfies \eqref{eqn:sqfree}
and $E/\Q$ is an elliptic curve such that
$\overline{\rho}_{F,n} \sim \overline{\rho}_{E,n}$.
In particular, for any prime $\ell \nmid N^\prime$,
\[
	\begin{cases}
		a_\ell(E) \equiv a_\ell(F) \mod{n} & \text{if $\ell \nmid N$}\\
		a_\ell(E) \equiv \pm (\ell+1) \mod{n} & \text{if $\ell \mid N$}.
	\end{cases}
\]

\subsection{The Method of Kraus}
\begin{lem}\label{lem:ratio}
Let $c^\prime = \pm c$ with the sign chosen so that $c^\prime \equiv 1 \mod{4}$.
Let
\[
\gamma=u+v \sqrt{-d} \; \; \mbox{ where } \; \; 
(u,v) \; = \; 
\begin{cases}
(1/8,3/8) & d=7\\
(7/8,1/8) & d=15\\
(3/8,1/8) & d=55\\
(5/16,-1/16) & d=231.
\end{cases}
\]
Choose $\epsilon_n \in \{1,-1\}$ to satisfy $n \equiv \epsilon_n \mod{3}$.
Then there is some $\delta \in M^*$ such that
\begin{equation}\label{eqn:fracprincipal}
\frac{x +c^\prime  \sqrt{-d}}{x -c^\prime \sqrt{-d}}  \; = \;
\begin{cases}
 \gamma \cdot \delta^n &  \mbox{ if }  d=7,15,55\\
\gamma^{(2+\epsilon_n \cdot n)/3} \cdot \delta^n & \mbox{ if }  d=231.
\end{cases}
\end{equation}
Moreover, $\delta$ is supported only on prime ideals dividing $y$.
\end{lem}
\begin{proof}
From the proof of Lemma~\ref{lem:facteven}, and in particular \eqref{eqn:fpfa},
we have
\begin{equation}\label{eqn:idealratio}
\left( \frac{x^\prime +c \sqrt{-d}}{x^\prime -c \sqrt{-d}} \right) \cdot \OO_M \; = \; (\overline{\fP}/\fP)^{2} \cdot \fB^n
\end{equation}
with $\fB=\fA/\overline{\fA}$. 
Here $\fP$ is given by \eqref{eqn:fP}, and $\fA$ is an integral
ideal dividing $y$.
We observe that  $\fB$ is supported only on prime ideals dividing $y$. 
First let $d=7$, $15$ or $55$. In these cases
the fractional ideal $(\overline{\fP}/\fP)^2$ is principal,
and we have chosen $\gamma$ so that it is a generator. Since $n$ is a prime
not dividing the order of $\Cl(M)$ we have that $\fB$ is also principal.
Let $\delta \in M^*$ be a generator of $\fB$. 
Then
\[
\frac{x +c^\prime \sqrt{-d}}{x -c^\prime \sqrt{-d}} = \pm \gamma \cdot \delta^n,
\]
and we complete the proof for $d=7$, $15$ and $55$ by absorbing the $\pm$ sign into $\delta$.

\medskip

Suppose now that $d=231$. The class of the fractional ideal $\overline{\fP}/\fP$
has order $3$, and we have chosen $\gamma$ to be
a generator of $(\overline{\fP}/\fP)^3$. 
We may rewrite \eqref{eqn:idealratio} as
\[
\frac{x +c^\prime \sqrt{-d}}{x -c^\prime \sqrt{-d}} \; = \; (\overline{\fP}/\fP)^{2+\epsilon_n \cdot n} \cdot \fC^n
\]
where $\fC=\fB \cdot (\fP/\overline{\fP})^{\epsilon_n}$. Note that $3 \mid (2+\epsilon_n \cdot n)$ and hence
\[
(\overline{\fP}/\fP)^{2+\epsilon_n \cdot n} \; = \; \gamma^{(2+\epsilon_n \cdot n)/3} \cdot \OO_M.
\]
The ideal $\fC$ must be principal and hence we complete the proof by letting $\delta$
be a suitably chosen generator for $\fC$. We note that, in all cases, $\delta$
is supported only on primes of $\OO_M$ dividing $y$.
\end{proof}

\begin{lem}\label{lem:Kraus}
	Let $n \ge 13$ be a prime and $(E,d)$ be one of the 
	remaining $114$ pairs.
Let $q=kn+1$ be a prime.
Suppose that $(-d/q)=1$, and choose $a$ such that $a^2 \equiv -d \mod{q}$.
Let $g_0$ be a generator for $\F_q^*$ and $g=g_0^n$.
Let $(u,v)$ be as in the statement of Lemma~\ref{lem:ratio}.
If $d=7$, $15$ or $55$, then let 
\[
\Theta_q^\prime \; = \;  \left\{   (u+v a) \cdot g^i \; : \; i=0,1,\dotsc,k-1\right\} \subset \F_q.
\]
If $d=231$, then set 
\[
\Theta_q^\prime \; = \;  \left\{   (u+v a)^{(2+\epsilon_n \cdot n)/3} \cdot g^i \; : \; i=0,1,\dotsc,k-1\right\} \subset \F_q
\]
and, in all cases, let
\[
\Theta_q \; = \; \Theta_q^\prime \setminus \{0,1\}.
\]
Suppose the following two conditions hold:
\begin{enumerate}
\item[(i)] $a_q(E)^2 \not\equiv 4 \mod{n}$.
\item[(ii)] $a_q(E)^2 \not\equiv a_q(H_\theta)^2 \mod{n}$ for all $\theta \in \Theta_q$, where
\[
H_\theta \;\; : \;\; Y^2=X(X+1)(X+\theta).
\] 
\end{enumerate}
Then $\overline{\rho}_{F,n} \nsim \overline{\rho}_{E,n}$.
\end{lem}
\begin{proof}
We suppose that $\overline{\rho}_{F,n} \sim \overline{\rho}_{E,n}$ and derive a contradiction.
Since $n \ge 11$, we note that, in particular,  $q \not\in \{  2, 3, 5, 7, 11 \}$.
Suppose first that $q \mid y$. Then  $q+1 \equiv \pm a_q(E) \mod{n}$.
But $q+1=kn+2 \equiv 2 \mod{n}$ and hence $a_q(E)^2 \equiv 4 \mod{n}$,
contradicting hypothesis (i).
We may therefore suppose that $q \nmid y$. In particular $q$ is a prime
of good reduction for the Frey curve $F$, and also for the curve $E$,
whence $a_q(F) \equiv a_q(E) \mod{n}$.

Since $a^2 \equiv -d \mod{q}$, by the Dedekind-Kummer theorem, the prime
$q$ splits in $\OO_M$ as a product of two primes $q \OO_M = \fq \cdot \overline{\fq}$
where we choose
\begin{equation}\label{eqn:fq}
\fq \; = \;  q \OO_M + (a-\sqrt{-d}) \cdot \OO_M.
\end{equation}
In particular $a \equiv \sqrt{-d} \mod{\fq}$. 
Moreover, $\F_\fq =\F_q$.
%As $\F_{\fq} = \F_q$, we have $a_{\fq}(F) \equiv a_q(E) \mod{n}$.
Since $\fq \mid q$ and $q \nmid 2y$, it follows 
from \eqref{eqn:rewrite} that  $\fq \nmid (x\pm c^\prime \sqrt{-d})$.
We let $\theta \in \F_q^*$ satisfy
\begin{equation}\label{eqn:thetadef}
\theta \; \equiv \; \frac{x+c^\prime\sqrt{-d}}{x-c^\prime\sqrt{-d}} \mod{\fq}.
\end{equation}
We will contradict hypothesis (ii), and complete the proof, by showing that
$\theta \in \Theta_q$ and $a_q(F) = \pm a_q(H_\theta)$.
If $\theta \equiv 1 \mod{q}$ then $\fq \mid 2 c^\prime \sqrt{-d}$
giving that $q \mid 2 \cdot 3 \cdot 5 \cdot 7 \cdot 11$, which is impossible.
Therefore $\theta \not \equiv 1 \mod{q}$. 
Let $(u,v)$, $\gamma$ and $\delta$
be as in the statement of Lemma~\ref{lem:ratio}. 
Note that $\gamma$ is supported
only at the primes above $2$ and that $\delta$ is supported
at only the primes above $y$. Since $\fq  \nmid y$,
we may reduce $\gamma$ and $\delta$ modulo $\fq$. In particular,
\[
\gamma \equiv u+a v \mod{\fq}.
\]
Moreover, $\delta^n \mod{\fq}$ belongs to the subgroup of $\F_q^*$
generated by $g=g_0^{n}$ of order $k$. The fact that $\theta$
belongs to $\Theta_q^\prime$ (and therefore to $\Theta_q$)
follows from \eqref{eqn:thetadef} and \eqref{eqn:fracprincipal}.

It remains to show that $a_q(F) = \pm a_q(H_\theta)$.
The model for $F$ in \eqref{eqn:evenFrey} is isomorphic to the model
in \eqref{eqn:simplifiedModel}.
We note that the polynomial on the right hand-side of \eqref{eqn:simplifiedModel}
can be factored as
\begin{equation}\label{eqn:Freyfactor}
X\left(X \, + \, \frac{x+c^\prime\sqrt{-d}}{2} \right) \left(X \, + \, \frac{x-c^\prime\sqrt{-d}}{2} \right).
\end{equation}
Thus, $F \mod \fq$ is a quadratic twist of $H_\theta$, whence
$$
a_q(F) = a_{\fq}(F)=\pm a_{\fq}(H_\theta)=\pm a_q(H_\theta),
$$
completing the proof.
\end{proof}

\noindent \textbf{Remark.}
We know by Dirichlet's theorem that the natural
density  
of primes $q$ satisfying the conditions $q=kn+1$
and $(-d/q)=1$ is $1/2n$. We now give a heuristic
estimate for the probability of succeeding
to show that $\overline{\rho}_{F,n} \nsim \overline{\rho}_{E,n}$
using a single $q=kn+1$ that satisfies $(-d/q)=1$.
The set $\Theta_q^\prime$ has size $k$, and so $\Theta_q$
has size close to $k$.
For a given $\theta \in \Theta_q$, we expect the probability
that $a_q(E)^2 \not\equiv a_q(H_\theta)^2 \mod{n}$
to be roughly $(1-2/n)$. Thus the probability of the criterion
succeeding is around $(1-2/n)^k$. In particular, if $k$ is large
compared to $n/2$ then we expect failure, but if $k$ is small
compared to $n/2$ then we expect success. 
Moreover, if we fail with one particular value of $q$,
we are likely to fail with larger values of $q$
(which correspond to larger values of $k$).

However,
this heuristic is likely to be inaccurate when
$\sqrt{q}$ is small compared to $n$, since
$a_q(E)$ and $a_q(H_\theta)$ both belong
to the Hasse interval $[-2\sqrt{q},2\sqrt{q}]$,
and the probability of the criterion succeeding is
around $(1-1/\sqrt{q})^k$.
\begin{comment}
When $\sqrt{q}$ is small compared to $n$,
then $a_q(E)^2 \equiv a_q(H_\theta)^2 \mod{n}$
is equivalent to $a_q(E) = \pm a_q(H_\theta)$.
Let us suppose additionally that $(\Delta(E)/q)=1$
where $\Delta(E)$ is the discriminant of $E$.
Since $E$ already has non-trivial $2$-torsion,
this additional assumption implies that $E/\F_q$
has full $2$-torsion. 
Thus $a_q(E)$ and $a_q(H_\theta)$ are integers
belonging to $[-2\sqrt{q},2\sqrt{q}]$
and satisfying $q+1-a_q(E) \equiv q+1-a_q(H_\theta) \equiv 0 \mod{4}$,
i.e. $a_q(E) \equiv a_q(H_\theta) \mod{4}$.
Thus we expect that the probability that $a_q(E)=\pm a_q(H_\theta)$
to be 
Since $E$ and $H_\theta$ are elliptic curves
with non-trivial $2$-torsion, the probability
that $a_q(E)=\pm a_q(H_\theta)$ is 
$1/\sqrt{q}$. Hence the probability
of the criterion succeeding is
$(1-1/\sqrt{q})^k$.
In this case we expect success if 
$k$ is small compared to $\sqrt{q}$
and failure if $k$ is large compared to $\sqrt{q}$.
\end{comment}

\begin{table}
\begin{tabular}{|c|c|c|c|}
	\hline
	$d$ & $N(d)$ & Number of primes $13 \le n < N(d)$ & Number of pairs $(E,d)$\\
	\hline\hline
	$7$ & $6 \times 10^8$ & $31324698$ & $39$\\
	\hline
	$15$ & $4 \times 10^8$ & $21336321$ & $28$\\
	\hline
	$55$ & $5 \times 10^8$ & $26355862$ & $27$\\
	\hline
	$231$ & $1.2 \times 10^9$ & $60454700$ & $20$\\
	\hline
\end{tabular}
	\caption{The upper bounds $N(d)$
	are as in \eqref{n-upper-bound}.
	The table records the number
	of primes in the interval
	$13 \le n < N(d)$
	and the number of pairs $(E,d)$.}\label{table:counts}
\end{table}

\bigskip

We are working towards proving Proposition~\ref{prop:yeven}.
Recall that we have $114$ remaining pairs $(E,d)$
with $E/\Q$ an elliptic curve and $d \in \{7,15,55,231\}$; these are 
distributed among the values of $d$ according to Table~\ref{table:counts}.
The table also records the upper bounds $N(d)$ of Proposition~\ref{prop:yeven}.
%There are $21336321$ primes $n$ in the 
%interval $13 \le n < 4 \times 10^8$, $26355862$ in the 
%interval $13 \le n < 5 \times 10^8$, $31324698$ in the 
%interval $13 \le n < 6 \times 10^8$, and $60454700$ in the 
%interval $13 \le n < 1.2 \times 10^9$.
%Recall that we have $114$ remaining pairs $(E,d)$
%with $E/\Q$ an elliptic curve and $d \in \{7,15,55,231\}$; these are distributed with $39, 28, 27$ and $20$  curves corresponding to $d=7, 15, 55$ and $231$, respectively.
We wrote a \texttt{Magma} script that applied the criterion of Lemma~\ref{lem:Kraus}
to the
$$
39 \cdot 31324698+28 \cdot 21336321+27 \cdot 26355862+20 \cdot 60454700
= 3739782484 \approx 3.7 \times 10^9
%28 \cdot 11078932 + 27 \cdot 13679313 +  39 \cdot 16252320 +  20 \cdot 31324698 = 1939885987 \approx 1.9 \times 10^9
$$
triples $(E,d,n)$. For each such triple, the script searches
for a prime $q=kn+1$ with $k <10^3$
such that the hypotheses of Lemma~\ref{lem:Kraus}
are satisfied. This computation took around
$29000$ hours, but was in fact distributed
over $64$ processors, and finished
in around $20$ days. For all but $1230$
of the $3739782484$ triples $(E,d,n)$ the 
script found some $q$ satisfying the hypotheses
of Lemma~\ref{lem:Kraus}. We are therefore
reduced to considering the remaining $1230$
triples $(E,d,n)$. While these are somewhat too numerous  to record
here, we note that the largest value of $n$ appearing
in any of these triples is $n=1861$ 
and this corresponds to $E$ being the
elliptic curve with Cremona label \texttt{210A1}
and $d=15$.

\subsection{A refined sieve}
Our adaptation of the method of Kraus (Lemma~\ref{lem:Kraus})
makes use of one auxiliary prime $q$ satisfying $q=kn+1$
and $(-k/q)=1$. To treat the remaining $1230$
triples $(E,d,n)$, we will use
a refined sieve that combines information
from several such primes $q$.

\begin{lem}\label{lem:sieve}
	Let $(E,d,n)$ be one of the 
	remaining $1230$ triples.
Let $q=kn+1$ be a prime.
Suppose that $(-d/q)=1$ and choose $a$ such that $a^2 \equiv -d \mod{q}$.
Let $c^\prime$, $h$, $(r,s)$, $\kappa_n$ be as in Lemma~\ref{lem:facteven},  $m=(2+ \kappa_n \cdot n)/h \in \Z$, and set
\[
\rho_1=(r+sa)^m \; \; \mbox{ and } \; \;  \rho_2=(r-sa)^m.
\]
Let $g_0$ be a generator for $\F_q^*$ and set $g=g_0^n$.
Further, let us define
\[
\Upsilon_q^{\prime\prime} \; = \;  \left\{   \left( \rho_1 \cdot g^i, \rho_2 \cdot g^j \right) \; : \; i=0,1,\dotsc,k-1,~j=0,1\right\} \subset \F_q \times \F_q,
\]
\[
\Upsilon_q^{\prime} \; = \;  \left\{   ( \theta_1,\theta_2) \in \Upsilon_q^{\prime\prime} \; : \; \theta_1 \theta_2 (\theta_1-\theta_2) \ne 0\right\}
\]
and
\[
\Upsilon_q \; = \;  \left\{   ( \theta_1,\theta_2) \in \Upsilon_q^{\prime} \; : \; a_q(H_{\theta_1,\theta_2})  \equiv a_q(E) \mod{n}\right\},
\]
where $H_{\theta_1,\theta_2}/\F_q$ is the elliptic curve
\[
H_{\theta_1,\theta_2} \; : \; Y^2 = X(X+\theta_1)(X+\theta_2).
\]
Write
\[
\Phi_q^\prime \; = \; \left\{ (\theta_1-\theta_2)/a \cdot (\F_q^*)^{2n} \; : \; (\theta_1,\theta_2) \in \Upsilon_q \right\} \subset \F_q^*/(\F_q^*)^{2n}
\]
and
\[
\Phi_q \; = \; 
\begin{cases}
\Phi_q^\prime \cup \{ (\omega/a) \cdot (\F_q^*)^{2n} \; : \; \omega \in \{ \rho_1,~\rho_1 g,~-\rho_2,~-\rho_2 g \} \} & \text{if $a_q(E)^2 \equiv 4 \mod{n}$}\\
\hskip22ex \Phi_q^\prime & \text{otherwise}.
\end{cases}
\]
If $\overline{\rho}_{F,n} \sim \overline{\rho}_{E,n}$, then necessarily 
\begin{equation}\label{eqn:cprimeclass}
c^\prime \cdot (\F_q^*)^{2n} \;   \in \;  \Phi_q.
\end{equation}
\end{lem}
\begin{proof}
Let $M=\Q(\sqrt{-d})$ and  $\fq \mid q$ be the  prime ideal of $\OO_M$ given by \eqref{eqn:fq}, so that $\OO_M/\fq=\F_q$ and $\sqrt{-d} \equiv a \mod{\fq}$.
Let $\mu$ be as in Lemma~\ref{lem:facteven}. 
From \eqref{eqn:principal1} and its conjugate, we have
\begin{equation}\label{eqn:rhoimu}
\frac{x+c^\prime \sqrt{-d}}{2}  \equiv \rho_1 \cdot \mu^n \mod{\fq}  \; \; \mbox{ and } \; \; 
\frac{x-c^\prime \sqrt{-d}}{2}  \equiv \rho_2 \cdot \overline{\mu}^n \mod{\fq}.
\end{equation}
Suppose first that $q \nmid y$. Thus both $F$ and $E$ have good reduction at $q$,
and so $a_q(F) \equiv a_q(E) \mod{n}$. It follows from \eqref{eqn:rewrite} that
$\fq \nmid ((x\pm c^\prime\sqrt{-d})/2)$ and that $\fq \nmid \mu$, $\overline{\mu}$.
Recall that $g=g_0^n$ where $g_0$ is a generator for $\F_q^*$;
	in particular, $g$ is a non-square, it generates $(\F_q^*)^n$, and has
	order $k$.
We note that the class of $\overline{\mu}^n$ modulo $\fq$ is either in $(\F_q^*)^{2n}$ or
in $g \cdot (\F_q^*)^{2n}$. Hence there is some
$\phi \in (\F_q^*)^{2n}$
and some $0 \le j \le 1$ such that 
\[
\frac{x-c^\prime \sqrt{-d}}{2}  \equiv \rho_2 \cdot g^j \cdot \phi \mod{\fq}.
\]
Now the class of $\mu^n/\phi$ modulo $\fq$ belongs to $(\F_q^*)^n$ and so
is equal to $g^i$ for some $0 \le i \le k-1$. We note that
\[
\frac{x+c^\prime \sqrt{-d}}{2}  \equiv \rho_1 \cdot g^i \cdot \phi \mod{\fq}.
\]
Hence
\[
\left( \frac{x+c^\prime \sqrt{-d}}{2} \, , \,  \frac{x-c^\prime \sqrt{-d}}{2} \right) \; \equiv\; (\theta_1 \cdot \phi , \theta_2 \cdot \phi) \mod{\fq}
\]
where $(\theta_1,\theta_2) \in \Upsilon_q^{\prime\prime}$. Since $\fq \nmid ((x\pm c^\prime\sqrt{-d})/2)$, we see that $\theta_1 \theta_2 \ne 0$.
Moreover, $\theta_1 -\theta_2 = c^\prime \sqrt{-d}/\phi \in \F_q^*$. Thus $(\theta_1,\theta_2) \in \Upsilon_q^\prime$.
Now recall that the model for the Frey curve $F$ in \eqref{eqn:evenFrey} is isomorphic
to the model given in \eqref{eqn:simplifiedModel}. The polynomial on the right hand-side
of the latter model factors as in \eqref{eqn:Freyfactor}. Thus $F/\F_\fq$
is isomorphic to the elliptic curve
\[
Y^2= X(X+\theta_1 \phi) (X+\theta_2 \phi).
\]
As $\phi$ is a square in $\F_q$, we see that this elliptic curve
is in turn isomorphic to the elliptic curve $H_{\theta_1,\theta_2}$.
Then $a_q(F)=a_\fq(F)=a_\fq(H_{\theta_1,\theta_2})=a_q(H_{\theta_1,\theta_2})$.
Since $a_q(E) \equiv a_q(F) \mod{n}$, it follows that
$(\theta_1,\theta_2) \in \Upsilon_q$. Moreover,
\[
c^\prime  \; = \; \frac{1}{\sqrt{-d}} \cdot \left( \frac{x+c^\prime \sqrt{-d}}{2} \; - \; \frac{x-c^\prime \sqrt{-d}}{2} \right)  \; \equiv \; \frac{\theta_1-\theta_2}{a} \cdot \phi \mod{\fq}.
\]
Since $\phi \in (\F_q^*)^{2n}$, this proves \eqref{eqn:cprimeclass}.

So far we have considered only the case $q \nmid y$. We know that if $q \mid y$, then 
$$
a_q(E) \equiv \pm (q+1) \equiv \pm 2 \mod{n}.
$$
Thus if $a_q(E)^n \not \equiv 4 \mod{n}$, then $q \nmid y$ and the proof is complete. Suppose $a_q(E)^2 \equiv 4 \mod{n}$
and that $q \mid y$. In particular, either $\fq \mid \mu$ or $\fq \mid \overline{\mu}$, but not both (by the coprimality of the 
factors on the right hand-side of \eqref{eqn:rewrite}).
Suppose $\fq \mid \overline{\mu}$. Then $x \equiv c^\prime \sqrt{d} \mod{\fq}$ and so
from \eqref{eqn:rhoimu} we have
\[
c^\prime \; \equiv \; \frac{\rho_1}{\sqrt{-d}}  \cdot \mu^n \; \equiv\;  \frac{\rho_1}{a} \cdot \mu^n \mod{\fq}.
\]
However, the class of $\mu^n$ modulo $\fq$ belongs to either $(\F_q^*)^{2n}$ or $g \cdot (\F_q^*)^{2n}$,
establishing \eqref{eqn:cprimeclass}. The case $\fq \mid \mu$ is similar. This completes the proof.
\end{proof}
\begin{lem}\label{lem:sieve2}
	Let $(E,d,n)$ be one of the 
	remaining $1230$ triples.
Let $q=kn+1$ be a prime.
Suppose that $(-d/q)=-1$. Let $M=\Q(\sqrt{-d})$
and let $\fq=q \OO_M$. Write $\F_{\fq}=\OO_M/\fq \cong \F_{q^2}$.
Let $c^\prime$, $h$, $(r,s)$, $\kappa_n$ be as in Lemma~\ref{lem:facteven}, and set $m=(2+ \kappa_n \cdot n)/h \in \Z$. Define
$\rho_1=(r+sa)^m$, choose
$g_0$ to  be a generator for $\F_\fq^*$, and set $g=g_0^n$.
Define
\[
\Upsilon_q^{\prime\prime} \; = \;  \left\{   \rho_1 \cdot g^i \; : \; i=0,1,\dotsc,2q+1\right\} \subset \F_\fq^*, 
\]
\[
\Upsilon_q^{\prime} \; = \;  \left\{   \theta \in \Upsilon_q^{\prime\prime} \; : \; \theta \ne \theta^q\right\}
\]
and
\[
\Upsilon_q \; = \;  \left\{   \theta \in \Upsilon_q^{\prime} \; : \; a_q(H_{\theta})  \equiv a_q(E) \mod{n}\right\},
\]
where $H_{\theta}/\F_q$ is the elliptic curve
\[
H_{\theta} \; : \; Y^2 = X(X+\theta)(X+\theta^q).
\]
Let
\[
\Phi_q \; = \; \left\{ (\theta-\theta^q)/\sqrt{-d} \cdot (\F_q^*)^{2n} \; : \; \theta \in \Upsilon_q \right\} \subset \F_q^*/(\F_q^*)^{2n}.
\]
If $\overline{\rho}_{F,n} \sim \overline{\rho}_{E,n}$ then necessarily 
 \eqref{eqn:cprimeclass} holds.
\end{lem}
\begin{proof}
We note that in $\F_\fq$ Galois conjugation agrees with the action
of Frobenius. Thus if $\alpha \in \OO_M$ and $\overline{\alpha}$ denotes its conjugate, then
 $\overline{\alpha} \equiv \alpha^q \pmod{\fq}$.

Since $(-d/q)=-1$ and $x^2+c^2d =y^n$ we observe that $q \nmid y$. Thus $F$ and $E$
both have good reduction at $q$, and so $a_q(F) \equiv a_q(E) \pmod{n}$.
Let $\mu$ be as in Lemma~\ref{lem:facteven}. 
Thus $\fq \nmid \mu$, $\overline{\mu}$. 
Recall that $g=g_0^n$ where $g_0$ is a generator for $\F_\fq^*$, whence $\mu^n \equiv g^j$ for some integer $j$. 
From \eqref{eqn:principal1},
\[
\frac{x+c^\prime \sqrt{-d}}{2} \equiv \rho_1 \cdot g^j \mod{\fq} \; \;  \mbox{ and } \; \; 
\frac{x-c^\prime \sqrt{-d}}{2} \equiv (\rho_1 \cdot g^j)^q \mod{\fq}.
\]
Write $j=i+(2q+2)t$,  where $i \in \{0,1,\dotsc,2q+1\}$ and $t$ is an integer.
We note that
\[
g^{2q+2}=(g_0^{q+1})^{2n}.
\]
Moreover, $g_0^{q+1}=g_0 g_0^{q} \in \F_q^*$. Thus there is some $\theta \in \Upsilon_q^{\prime\prime}$
and some $\phi \in (\F_q^*)^{2n}$ such that
\[
\frac{x+c^\prime \sqrt{-d}}{2} \equiv \theta \cdot \phi \mod{\fq} \; \; \mbox{ and } \; \;  \frac{x-c^\prime\sqrt{-d}}{2} \equiv \theta^q \cdot \phi \mod{\fq}.
\] 
Since $\fq \nmid c^\prime \sqrt{-d}$, we see that $\theta \ne \theta^q$ and so $\theta \in \Upsilon_q^\prime$.
We note that the model for $F$ in \eqref{eqn:simplifiedModel} can, over $\F_q$, be written as
\[
Y^2 = X(X^2+ \phi \cdot (\theta+\theta^q) X+ \phi \cdot (\theta \theta^q)),
\]
where the coefficients are fixed by Frobenius and so do indeed belong to $\F_q$.
This model is a twist by $\phi$ of $H_\theta$. As $\phi$ is a square in $\F_q^*$, 
we have $a_q(H_\theta)=a_q(F) \equiv a_q(E) \pmod{n}$. Thus $\theta \in \Upsilon_q$.
Finally, 
\[
c^\prime  \; = \; \frac{1}{\sqrt{-d}} \cdot \left( \frac{x+c^\prime \sqrt{-d}}{2} \; - \; \frac{x-c^\prime \sqrt{-d}}{2} \right)  \; \equiv \; \frac{\theta-\theta^q}{\sqrt{-d}} \cdot \phi \mod{\fq}.
\]
Since $\phi \in (\F_q^*)^{2n}$,  this proves \eqref{eqn:cprimeclass}.
\end{proof}

\begin{lem}\label{lem:refined}
	Let $(E,d,n)$ be one of the 
	remaining $1230$ triples.
	Let $q_1,q_2,\dotsc,q_r$ be primes satisfying
		$q_i \equiv 1 \mod{n}$.
Let
\[
	\psi_q \; : \; (\Z/2n \Z)^4 \rightarrow \F_q^*/(\F_q^*)^{2n}, \qquad \psi_q(x_1,x_2,x_3,x_4)=(-3)^{x_1} 5^{x_2} (-7)^{x_3} (-11)^{x_4} \cdot (\F_q^*)^{2n}.
\]
	If $(-d/q)=1$, let $\Phi_{q_i}$ be as in Lemma~\ref{lem:sieve}
and if $(-d/q)=-1$, let $\Phi_{q_i}$ be as in Lemma~\ref{lem:sieve2}.
	Suppose
	\[
		\bigcap_{i=1}^{r} \psi_{q_i}^{-1} (\Phi_{q_i}) \; = \; \emptyset.
	\]
	Then $\overline{\rho}_{F,n} \nsim \overline{\rho}_{E,n}$.
\end{lem}
\begin{proof}
%Lemma~\ref{lem:sieve} gives information about the class of $c^\prime$ in $(\F_q)^*/(\F_q^*)^n$.
Recall, from \eqref{eqn:yeven} and \eqref{eqn:sqfree}, that
\[
	c=3^{\beta_3} 5^{\beta_5} 7^{\beta_7} 11^{\beta_{11}}.
\]
Thus $c \equiv (-1)^{\beta_3+\beta_7+\beta_{11}} \mod{4}$ and hence, since we choose $c^\prime = \pm c$
so that $c^\prime \equiv 1 \mod{4}$,
\[
	c^\prime = (-1)^{\beta_3+\beta_7+\beta_{11}} \cdot 3^{\beta_3} 5^{\beta_5} 7^{\beta_7} 11^{\beta_{11}}
	= (-3)^{\beta_3} 5^{\beta_5} (-7)^{\beta_7} (-11)^{\beta_{11}}.
\]
%As before, let $q=kn+1$ be a prime satisfying $(-d/q)=1$, and choose some $a \in \F_q$
%so that $a^2 \equiv -d \mod{q}$. Let $\fq$ be given by \eqref{eqn:fq}.
Suppose $\overline{\rho}_{F,n} \sim \overline{\rho}_{E,n}$. 
Thus
\[
	\psi_q(\beta_3,\beta_5,\beta_7,\beta_{11})= c^\prime \cdot (\F_{q_i}^*)^{2n} \in \Phi_{q_i}
\]
	by \eqref{eqn:cprimeclass}. Therefore
	\[
		((\beta_3,\beta_5,\beta_7,\beta_{11}) \bmod{2n}) \; \in \; 
		\bigcap_{i=1}^{r} \psi_{q_i}^{-1} (\Phi_{q_i})
	\]
	giving a contradiction.
\end{proof}

We wrote a \texttt{Magma} script which for each of the $1230$ remaining triples $(E,d,n)$  
recursively computes the intersections
\[
\psi_{q_1}^{-1}(\Phi_{q_1}), \quad
\bigcap_{i=1}^2 \psi_{q_i}^{-1} (\Phi_{q_i}), \quad
\bigcap_{i=1}^3 \psi_{q_i}^{-1} (\Phi_{q_i}), \dots
\]
where the $q_i$ are primes $\equiv 1 \pmod{n}$. It stops when the intersection is empty,
or when we have used $200$ primes $q_i$, whichever comes first. If the intersection is empty,
then we know from Lemma~\ref{lem:refined} that $\overline{\rho}_{F,n} \nsim \overline{\rho}_{E,n}$
and we may eliminate the particular triple $(E,d,n)$ from further consideration. We reached an empty
intersection in $1224$ cases. Table~\ref{table:bad} gives
the details for the six triples $(E,d,n)$ where the intersection is non-empty.

\begin{table}
\caption{This table gives the six triples
$(E,d,n)$ such that the intersection
$\bigcap_{i=1}^{200} \psi_{q_i}^{-1} (\Phi_{q_i})$ is non-empty.
Here the elliptic curve $E$ is given in the first column in Cremona
notation.
We note that $n=13$ for all six triples. Therefore the intersection
given in the last column is a subset of $(\Z/26\Z)^4$.}
\label{table:bad}
\begin{tabular}{|c|c|c|c|c|}
\hline
Elliptic Curve & $d$ & $n$ & $\displaystyle \bigcap_{i=1}^{200} \psi_{q_i}^{-1} (\Phi_{q_i})$ \\
\hline\hline
\texttt{462b1} & $231$ & $13$ &
    $\{ \; ( 7, 2, 19, 3 ),\;  ( 9, 1, 24, 9 ) \; \}$
\\
\hline
\texttt{462f1} & $231$ & $13$ & $ \{ \; ( 0, 15, 25, 13 ), \; ( 15, 18, 5, 0 )\; \}$
\\
\hline
\texttt{2310j1} &  $231$ & $13$ & $ \{ \; ( 11, 6, 6, 18 ),\; ( 24, 19, 19, 5 )\; \}$\\
\hline
\texttt{2310l1} & $231$ & $13$ & $\{ \; ( 10, 5, 22, 8 )\; \}$\\
\hline
\texttt{2310m1} & $231$ & $13$ & $\{ ( 5, 14, 11, 21 ),\; ( 7, 21, 19, 19 )\; \}$\\
\hline
\texttt{2310o1} & $15$ & $13$ & $\{\; ( 1, 0, 1, 1 ) \; \}$\\
\hline
\end{tabular}
\end{table}

\subsection{Proof of Proposition~\ref{prop:yeven}}
We now complete the proof of Proposition~\ref{prop:yeven}.
To summarise, Lemma~\ref{lem:TM} showed that the only solutions 
to \eqref{eqn:yeven} with exponent $n \in \{ 5, 7, 11 \}$
are the ones given in the statement of Proposition~\ref{prop:yeven}.
In view of the results of this section, it only remains to consider
the six triples $(E,d,n)$ given in Table~\ref{table:bad}.
To eliminate further cases, we make use of the following result of Halberstadt and Kraus \cite[Lemme 1.6]{HK}.
\begin{thm}[Halberstadt and Kraus]\label{thm:HK}
Let $E_1$ and $E_2$ be elliptic curves over $\mathbb{Q}$ 
and write $\Delta_j$ for the minimal discriminant of $E_j$.
Let $n \ge 5$ be a prime such that
$\overline{\rho}_{E_1,n} \sim \overline{\rho}_{E_2,n}$.
Let $q_1$, $q_2 \ne n$ be distinct primes of multiplicative reduction
for both elliptic curves such 
that $\ord_{q_i}(\Delta_j) \not \equiv 0 \pmod{n}$ for $i,j \in \{ 1,2 \}$. Then
\[
\frac{\ord_{q_1}(\Delta_1) \cdot \ord_{q_2}(\Delta_1)}{\ord_{q_1}(\Delta_2) \cdot \ord_{q_2}(\Delta_2)}
\]
is congruent to a square modulo $n$.
\end{thm}

%\begin{proof}[Proof of Proposition~\ref{prop:even}]
We shall use Theorem~\ref{thm:HK} and Lemma~\ref{lem:condequiv} to eliminate the first five of the six outstanding
triples $(E,d,n)$ given in Table~\ref{table:bad}. 
In all these cases $n=13$.
We know from the proof of Lemma~\ref{lem:refined} that
$(\beta_3,\beta_5,\beta_7,\beta_{11}) \mod{26}$ belongs to the intersection in the last column
of Table 1.

Consider the first triple, corresponding to the first row of the table. 
The $\beta_5 \equiv 1$ or $2 \pmod{26}$. But $\beta_5=\ord_5(c)$. 
Thus
$2\ord_5(c)+\ord_5(d)\equiv 2\beta_5+\ord_5(231) \equiv 2$ or $4 \pmod{13}$ and so by Lemma~\ref{lem:condequiv},
$5$ must divide the conductor of $E$ which is $462$ giving a contradiction.
The same argument eliminates the second triple.

Next we consider the third triple. Here $\beta_7 \equiv 6$ or $19 \pmod{26}$,
and so $\ord_7(c) \equiv \beta_7 \equiv 6 \pmod{13}$. Then
$2\ord_7(c)+\ord_7(d) \equiv 2 \beta_7+\ord_7(231) \equiv 0 \pmod{13}$. By Lemma~\ref{lem:condequiv},
$7$ does not divide the conductor of $E$ which is $2310$, again a contradiction.

We next consider the fourth triple. Here the elliptic curve $E$ with Cremona reference \texttt{2310l1}
has minimal discriminant
\[
\Delta_E= 2^4 \times 3^{12} \times 5^3 \times 7 \times 11.
 %<2, 4>, <3, 12>, <5, 3>, <7, 1>, <11, 1>
\]
We apply Theorem~\ref{thm:HK} with $E_1=F$, $E_2=E$, $q_1=2$ and $q_2=3$.
From the proof of Lemma~\ref{lem:condequiv} we have 
\[
\ord_2(\Delta_F) \equiv -12 \equiv 1 \pmod{13}, \qquad \ord_3(\Delta_F)=2\beta_3+\ord_3(231) \equiv2 \times 10+1 \equiv 8 \pmod{13}.
\]
Hence
\[
\frac{\ord_2(\Delta_F) \cdot \ord_3(\Delta_F)}{\ord_2(\Delta_E) \cdot \ord_3(\Delta_E)}
\equiv \frac{1 \times 8}{4 \times 12} \equiv 11 \pmod{13}
\]
which is a non-square modulo $13$, contradicting Theorem~\ref{thm:HK}.

Next we consider the fifth triple. Here there are two possibilities 
for $(\beta_3,\beta_5,\beta_7,\beta_{11})$. In the second possibility
we have $\beta_{7} \equiv 19 \pmod{26}$ which leads to a contradiction
via Lemma~\ref{lem:condequiv}. We focus on the first possibility.
The minimal discriminant of the curve $E$ is
\[
\Delta_E=2^4 \times 3^8 \times 5 \times 7^{3} \times 11.
 %<2, 4>, <3, 8>, <5, 1>, <7, 3>, <11, 1> ]
\]
We obtain a contradiction by applying Theorem~\ref{thm:HK}
with $q_1=2$ and $q_2=3$.

We are left with the last triple, which we have been unable
to eliminate by appealing to Theorem~\ref{thm:HK}
or Lemma~\ref{lem:condequiv}, or by further sieving.
In fact, \eqref{eqn:yeven} has the solution
\begin{equation}\label{eqn:n13sol}
	8143^2+3^3 \cdot 5 \cdot 7^2 \cdot 11^2 = 4^{13}.
\end{equation}
Here $n=13$, $d=15$ and $c=3 \cdot 7 \cdot 11$. We note
that the vector of exponents for this value of $c$
is $(\beta_3,\beta_5,\beta_7,\beta_{11})=(1,0,1,1)$
which agrees with the prediction in the last
column of the table.
Moreover, letting $x=-8143 \equiv 1 \pmod{4}$,
and $y^n=4^{13}$ in the Frey curve $F$ gives
the elliptic curve \texttt{2310o1}.
To complete the proof,
we need to solve \eqref{eqn:yeven}
with $d=15$ and $n=13$. We do this by reducing this case to a Thue-Mahler equation
using the approach in the proof of Lemma~\ref{lem:TM}. 
After possibly changing the sign of $x$ so that $x \equiv 1 \pmod{4}$, we 
have that
\[
	\frac{x+c^\prime\sqrt{-15}}{2} =\left(\frac{1-\sqrt{-15}}{8}\right) \left(r \, + \, s \cdot \frac{(1+\sqrt{-15})}{2}\right)^{13}, 
\]
where $y=r^2+rs+4s^2$ for some integers $r$ and $s$. Equating imaginary parts leads to the conclusion that 
$$
F_{13}(r,s) = \sum_{i=0}^{13} a_i r^{13-i} s^i =  \pm  4 \cdot 3^{\beta_3} \cdot 5^{\beta_5} \, \cdot\,  7^{\beta_7} \cdot 11^{\beta_{11}},
$$
where
$$
\begin{array}{|cc|cc|cc|} \hline
i & a_i & i & a_i & i & a_i \\ \hline
0 & 1 & 5 & 36036 & 10 & 195624 \\
1 & 0 & 6 & -34320 & 11 & -95160\\
2 & -312 & 7 & -226512 & 12 &  -51428 \\
3 & -1144 & 8 & -66924 & 13 & 924. \\
4 & 8580 & 9 & 340340 &  &  \\ \hline
\end{array}
$$

%\begin{multline*}
%	r^{13} - 312 r^{11} s^2 - 1144 r^{10} s^3 + 8580 r^9 s^4 + 36036 r^8 s^5 -
 %   34320 r^7 s^6 - 226512 r^6 s^7\\
%	- 66924 r^5 s^8 + 340340 r^4 s^9 +
%    195624 r^3 s^{10} - 95160 r^2 s^{11} - 51428 r s^{12} + 924 s^{13} \; = \; \pm  4 \cdot 3^{\beta_3} \cdot 5^{\beta_5} \, \cdot\,  7^{\beta_7} \cdot 11^{\beta_{11}}.
%\end{multline*}
We solved this Thue-Mahler equation using the \texttt{Magma} package associated to the paper \cite{GKMS}.
The only solution is with
%[ 0, -1, 1, 0, 1, 1 ]
\[
r=0, \quad s=\pm 1, \quad \beta_3=1, \quad \beta_5=0, \quad \beta_7=1 \; \; \mbox{ and } \; \;   \beta_{11}=1.
\]
This corresponds to the identity \eqref{eqn:n13sol} and completes
the proof of Proposition~\ref{prop:yeven}.

% The computation took 1314 s = 22 min

\bigskip

\noindent \textbf{Remark.}
It is natural to ask if the case $n=13$ could 
have been dealt with entirely using
the Thue-Mahler approach, just as we 
did for $n \in \{ 5, 7, 11 \}$ in Lemma~\ref{lem:TM}.
The Thue-Mahler solver that we
are using can quickly deal with 
the Thue-Mahler equations associated
to the pairs $(d,n)=(7,13)$ and $(55,13)$.
However, the Thue-Mahler equation
for the pair $(d,n)=(231,13)$
appears to be somewhat beyond its
capabilities. The approach in \cite{GKMS}
reduces solving a Thue-Mahler equation to
 a certain number of
$S$-unit equations. By way of example, the Thue-Mahler
equation for the pair $(d,n)=(15,13)$
reduces to solving four $S$-unit equations.
The Thue-Mahler equation for the pair
$(d,n)=(231,13)$, on the other hand,  corresponds to
$2240$ $S$-unit equations.
This explains the effort we
invested into eliminating 
$(d,n)=(231,13)$
via sieving and appeal to
Theorem~\ref{thm:HK}
and Lemma~\ref{lem:condequiv}.

\begin{comment}
\texttt{462b1} & $231$ & $13$ &
    $\{[ 7, 2, 19, 3 ], [ 9, 1, 24, 9 ]\}$
& 
    $\{[ 7, 2, 19, 3 ],
    $[ 9, 1, 24, 9 ]$
\\
+++++++++++++++++++++++++++++++
462f1 231
{
    [ 0, 15, 25, 13 ],
    [ 15, 18, 5, 0 ]
}
{}
+++++++++++++++++++++++++++++++
2310j1 231
{
    [ 11, 6, 6, 18 ],
    [ 24, 19, 19, 5 ]
}
{}
+++++++++++++++++++++++++++++++
2310l1 231
{
    [ 10, 5, 22, 8 ]
}
{}
+++++++++++++++++++++++++++++++
2310m1 231
{
    [ 5, 14, 11, 21 ],
    [ 7, 21, 19, 19 ]
}
{
    [ 7, 21, 19, 19 ]
}
+++++++++++++++++++++++++++++++
2310o1 15
{
    [ 1, 0, 1, 1 ]
}
{
    [ 1, 0, 1, 1 ]
}
+++++++++++++++++++++++++++++++
\end{comment}

%---------------------------------------------------------------------------------
\section{Equation (\ref{eq-main}) with $y$ even : large exponents} \label{sec:large}
%-----------------------------------------------------------------------------------

From the results of the preceding sections, it remains to solve equation (\ref{eq-main}) with $y$ even and exponent $n$ prime and
 \begin{equation} \label{starter!}
n > N(d),
\end{equation}
where $N(d)$ is as defined in (\ref{n-upper-bound}). We will accomplish this through (quite careful) application of bounds for linear forms in logarithms.
 
\subsection{Upper bounds for $n$ : linear forms in logarithms, complex and $q$-adic}
Our first order of business will be to produce an upper bound for the exponent $n$; initially it will be somewhat larger than $N(d)$. To this end, as it transpires, it will prove useful to have at our disposal a lower bound upon $y$.
From the discussion following Lemma \ref{lem:expbound}, we have that $\overline{\rho}_{F,n} \sim \overline{\rho}_{E,n}$ for $E/\mathbb{Q}$ with nontrivial rational $2$-torsion.

%\begin{lem} \label{bounders}
%We have that
%$\overline{\rho}_{F,n} \sim \overline{\rho}_{E,n}$
%where $E$
%is an elliptic curve over $\mathbb{Q}$ with conductor 
%$$
%N \in \{ 14, 30, 42, 70, 154, 210, 330, 462, 770, 2310  \}
%$$
%and nontrivial rational $2$-torsion, unless possibly $n \in \{ 7, 13 \}$ and $N \in \{ 770, 2310 \}$, or $n=11$ and $N \in \{ 110, 154, 770, 2310 \}$. 
%\end{lem}

%\begin{proof}
%This follows from an application of Lemma \ref{lem:irratbound-first}. If we let
%$\cC_{f}$ denote the ideal generated by $\cC_{f,\ell}$
%where $\ell$ runs through the primes $11< \ell < 100$,
%it follows  that $n \mid \cB_f = \mbox{Norm}(\cC_{f})$.
%The newforms where $\cB_f=0$ correspond to elliptic curves with nontrivial rational $2$-torsion. For each $n \in \{ 7, 11, 13 \}$, there is at least one newform $f$ of dimension exceeding $1$ with $n \mid  \cB_f$.
%\end{proof}

To begin, we will need to treat the case where $y$ in equation (\ref{eq-main}) has no odd prime divisors. Suppose that we have a solution to equation (\ref{eqn:yeven})
with $y=2^\kappa$ for $\kappa$ a positive integer. For the time being, we will relax our assumptions upon $n$ and suppose only that $n \geq 7$ is prime. Then the Frey-Hellegouarch curve $F$ has nontrivial rational $2$-torsion and conductor 
$$
N=2 \cdot 3^{\delta_3} 5^{\delta_5}  7^{\delta_7} 11^{\delta_{11}} \; \mbox{ where } \; \delta_i \in \{ 0, 1 \},
$$
so that
$$
N \in \{ 14, 30, 42, 66,  70, 154, 210, 330, 462, 770, 2310 \},
$$
and minimal discriminant 
$$
-2^{2 \kappa n -12} 3^{\alpha_3}5^{\alpha_5}7^{\alpha_7}11^{\alpha_{11}}. 
$$
A quick check of Cremona's tables reveals that we find such curves with minimal discriminant negative and divisible by precisely  $2^{2 \kappa n -12}$, with $n \geq 7$ prime, only for $18$ isomorphism classes of curves, given, in Cremona's notation, by
$$
\begin{array}{l}
14a4, \, 210b5, \, 210e1, \, 210e6, \, 330c1, \, 330c6, \, 330e4, \, 462a1, \, 462d1, \, 462e1,\\
 462g3, \, 770a1, \, 770e1, \, 770g3, \, 2310d4, \, 2310n1, \, 2310n6, \, 2310o1. \\
 \end{array}
$$
Most of these have $2 \kappa n -12=2$ and so $\kappa=1$ and $n=7$. Since $P(2^7-x^2) > 11$ for $1 \leq x < 11$ odd, only the curve 14a4 with $ \Delta = - 2^{2}  \cdot 7$ corresponds to a solution, arising from the identity $11^2 + 7 = 2^7$. Four more curves have $2 \kappa n -12=16$ and so $\kappa=2$ and $n=7$. Corresponding identities are
$$
7^2+3^3 \cdot 5 \cdot 11^2 = 2^{14}, \; 47^2 + 3^4 \cdot 5^2 \cdot 11 = 2^{14}, \; \; 103^2+ 3 \cdot 5^2 \cdot 7 \cdot 11 = 2^{14}, \; \; 117^2 +  5 \cdot  7^2 \cdot 11 = 2^{14},
$$
arising from the curves 330c1, 210e1, 2310n1 and 770e1, with discriminants
$$
- 2^{16} \cdot 3^{3} \cdot 5 \cdot 11^2, \; \;  - 2^{16} \cdot 3^{4} \cdot 5^2 \cdot 7, \; \; - 2^{16} \cdot 3 \cdot 5^2 \cdot 7 \cdot 11 \; \mbox{ and } \; - 2^{16} \cdot 5 \cdot 7^2 \cdot 11,
$$
respectively. Neither 462d1 nor 462e1 lead to any solutions while 2310o1, with discriminant $- 2^{40} \cdot 3^3 \cdot 5 \cdot 7^2 \cdot 11^2$, corresponds to the identity
$$
8143^2+ 3^3 \cdot 5 \cdot 7^2 \cdot 11^2 = 2^{26}.
$$

We may thus suppose that $y$ is divisible by an odd prime factor, provided $n \geq 17$.

\begin{lem} \label{ybound}
If $n \geq 17$ and $y$ is even, we have 
$$
y > 4n -4 \sqrt{2n} + 2.
$$
\end{lem}
\begin{proof}
By our preceding remarks,  there necessarily exists an odd prime $p \mid y$. Since $\overline{\rho}_{F,n} \sim \overline{\rho}_{E,n}$ where $E/\mathbb{Q}$ has nontrivial rational $2$-torsion, 
the fact that $\gcd (x,y)=1$ thus allows us to conclude that 
$$
a_p(E) \equiv \pm (p+1) \mod{n}.
$$
From the Hasse-Weil bounds, we have that $a_p(E)$ is bounded in modulus by $2 \sqrt{p}$, so that, using the fact that $a_p(E)$ is even, 
$$
n <  \frac{1}{2} (\sqrt{p}+1)^2 \leq \frac{1}{2} (\sqrt{y/2}+1)^2.
$$
The desired inequality follows.
\end{proof}

As before, define $c$ and $d$ via (\ref{eqn:sqfree}), where, since $y$ is even, 
$d \in \{ 7, 15, 55, 231 \}$, and let $c^\prime = \pm c$ with the sign chosen so that $c^\prime \equiv 1 \mod{4}$.
To derive an upper bound upon $n$, we will begin by using (\ref{eqn:fracprincipal}) to find a ``small'' linear form in logarithms. Specifically, let us define
\begin{equation} \label{lambda!}
\Lambda = \log \left( \frac{x +c^\prime  \sqrt{-d}}{x -c^\prime \sqrt{-d}}   \right).
\end{equation}
We prove
\begin{lem} \label{smallform}
If we suppose that
\begin{equation} \label{ass}
y^n > 100 \, c^2d, 
\end{equation} 
then
$$
\log \left|  \Lambda \right| < 0.75 +\log c + \frac{1}{2} \log d - \frac{n}{2} \log y.
$$
\end{lem}
 
\begin{proof}
Assumption (\ref{ass}), together with, say, Lemma B.2 of Smart \cite{Smart}, implies that 
$$
\left| \Lambda \right| \leq -10 \log (9/10) \left| \frac{x +c^\prime  \sqrt{-d}}{x -c^\prime \sqrt{-d}}   -1 \right| = 20 \log (10/9)  \frac{c \sqrt{d}}{y^{n/2}},
$$
whence the lemma follows.
\end{proof}

To show that $\log |\Lambda|$ here is indeed small, we first require an upper bound upon the exponents $\alpha_q$ in equation (\ref{eqn:yeven}). From (\ref{eqn:fracprincipal}), we have that 
\begin{equation} \label{ake}
\frac{2 \cdot c^\prime \sqrt{-d}}{x-c^\prime \sqrt{-d}} = 
\begin{cases}
 \gamma \cdot \delta^n -1 &  \mbox{ if } d \in \{ 7,15,55 \} \\
\gamma^{(2+\epsilon_n \cdot n)/3} \cdot \delta^n -1 & \mbox{ if }  d=231.
\end{cases}
\end{equation}
For prime $q$, let $\overline{\mathbb{Q}_q}$ denote an algebraic closure of the $q$-adic field $\mathbb{Q}_q$, and define $\nu_q$ to be the unique extension to $\overline{\mathbb{Q}_q}$ of the standard $q$-adic valuation over $\mathbb{Q}_q$, normalized so that $\nu_q(q)=1$. For any algebraic number $\alpha$ of degree $d$ over $\mathbb{Q}$, we define the {\it absolute logarithmic height} of $\alpha$ via the formula
\begin{equation}\label{eqn:htdef}
h(\alpha)= \dfrac{1}{d} \left( \log \vert a_{0} \vert + \sum\limits_{i=1}^{d} \log \max \left( 1, \vert \alpha^{(i)}\vert \right) \right), 
\end{equation}
where $a_{0}$ is the leading coefficient of the minimal polynomial of $\alpha$ over $\mathbb{Z}$ and the $\alpha^{(i)}$ are the conjugates of $\alpha$ in $\mathbb{C}$.
Since $\gcd (x,q)=1$, it follows from (\ref{ake}) that, if we set
$$
\Lambda_1 =
\begin{cases}
 \delta^n -(1/\gamma) &  \mbox{ if } d \in \{ 7, 15, 55 \} \\
 \delta^n -(1/\gamma)^{(2+\epsilon_n \cdot n)/3} & \mbox{ if }  d=231,
\end{cases}
$$
then $\nu_q(\Lambda_1) \geq \alpha_q/2$, for $q\in \{ 3, 5, 7, 11 \}$.

To complement this with an upper bound for linear forms in $q$-adic logarithms, we will appeal to Th\'eor\`eme 4 of Bugeaud  and Laurent \cite{BuLa}, with, in the notation of that result, the choices $(\mu,\nu)=(10,5)$.
\begin{thm}[Bugeaud-Laurent] \label{qlog}
Let $q$ be a prime number and let $\alpha_1, \alpha_2$ denote algebraic numbers which are $q$-adic units. Let $f$ be the residual degree of the extension $\mathbb{Q}_q(\alpha_1,\alpha_2)/\mathbb{Q}_q$ and put $D=[\mathbb{Q}_q(\alpha_1,\alpha_2) : \mathbb{Q}_q]/f$. Let $b_1$ and $b_2$ be positive integers and put
$$
\Lambda_1 = \alpha_1^{b_1}-\alpha_2^{b_2}.
$$
Denote by $A_1 > 1$ and $A_2 > 1$ real numbers such that
$$
\log A_i \geq \max \left\{ h(\alpha_i), \frac{\log q}{D} \right\}, \; \; i \in \{ 1, 2 \},
$$
and put
$$
b^\prime = \frac{b_1}{D \log A_2} + \frac{b_2}{D \log A_1}.
$$
If $\alpha_1$ and $\alpha_2$ are multiplicatively independent, then we have the bound
$$
\nu_q(\Lambda_1) \leq \frac{24 q (q^f-1)}{(q-1)\log^4 (q)} \, D^4 \left( \max \left\{ \log b^\prime + \log \log q + 0.4, \frac{10 \log q}{D}, 5 \right\} \right)^2 \cdot \log A_1  \cdot \log A_2.
$$
\end{thm}

%Arguing as in \cite{BMS2}, 

%and we have
%$$
%h(\beta_q) = \frac{1}{2} \, h_q \log 2 \; \; \mbox{ and } \; \; h(\gamma) = \frac{1}{2} \, h_q \log y.
%$$

%Let $0 \le r <h$ be such that
%$\fP^r \cdot \fA$ is principal, whereby, from \eqref{eqn:pretm},
%we deduce that $\fP^{n-2-nr}$ is principal. It follows that $r$
%is the unique solution to the congruence $n r \equiv 2-n \mod{h}$
%satisfying $0 \le r < h$, and we now fix $r$ to be that value.
%Let $\beta$ be a generator for the fractional ideal
%$\fP^{n-2-nr}$ and let $\gamma \in \OO_M$
%be a generator for $\fP^r \cdot \fA$. Note that the units
%of $\OO_M$ are $\pm 1$. After possibly changing the sign of $\gamma$
%we may rewrite \eqref{eqn:pretm} as
%\begin{equation}\label{eqn:pretm2}
%\frac{x+q^k\sqrt{-q}}{2}  = \beta \cdot \gamma^n,
%\end{equation}
%where $\beta$ is known but $\gamma$ is unknown.

We will choose $q \in \{ 3, 5, 7, 11 \}$ and apply this result  with the following choices of parameters :
$$
f=1, \;  \; D=2, \; \alpha_1= \delta, \;  \; \alpha_2 = 1/\gamma, \;  b_1=n
$$
and
\begin{equation} \label{bee2}
 b_2=
 \begin{cases}
\; \; \; \; \;  \; \; \; 1 &  \mbox{ if } d \in \{ 7, 15, 55 \} \\
 (2+\epsilon_n \cdot n)/3 & \mbox{ if }  d=231.
\end{cases}
 \end{equation}
We have 
 $$
 h(1/\gamma) = 
  \begin{cases}
 \log 2 &  \mbox{ if } d \in \{ 7, 15, 55 \} \\
 \frac{3}{2} \log 2 & \mbox{ if }  d=231
\end{cases}
 $$
 and
 $$
 h(\delta) \leq \frac{1}{2} \log (y/2),
 $$
 and hence, from (\ref{starter!}) and Lemma \ref{ybound},  may choose
$$
 \log A_1= \frac{1}{2} \log (y/2) <  \frac{1}{2} \log y
 $$
and 
$$
\log A_2 =  
\begin{cases}
\log 2  & \mbox{if } d \in \{ 7, 15, 55 \} \mbox{ and }  q=3, \\
 \frac{3}{2} \log 2 & \mbox{if }  d=231 \mbox{ and }  q \in \{ 3, 5, 7 \}, \\
 \frac{1}{2} \log q & \mbox{otherwise}. \\
\end{cases}
$$
Once again appealing to (\ref{starter!}) and Lemma \ref{ybound}, we have, in all cases,  that  $b^\prime > 5 \log q$ and 
$$
b^\prime \leq  \frac{n}{ 2 \log 2 } + \frac{1}{ \log y} < 0.722 n.
$$
We thus have
$$
\log b^\prime + \log \log q + 0.4 < 1.05 \log n,
$$
whence, from Theorem \ref{qlog},

$$
\nu_q(\Lambda_1) < c(d,q) \cdot  1.05^2 \log^2 n  \log y,
$$
where 
$$
c(d,q)=
\begin{cases}
 \frac{576 \log 2}{\log^4 3}  & \mbox{if } d \in \{ 7, 15, 55 \} \mbox{ and }  q=3,  \\
 \frac{288  q \log 2}{\log^4 q}  &  \mbox{if }  d=231 \mbox{ and }  q \in \{ 3, 5, 7 \}, \\
\frac{96 q}{\log^3 q}  & \mbox{otherwise}. \\
 \end{cases}
 $$
It follows that
\begin{equation} \label{whopper}
\sum_{q \in \{ 3, 5, 7, 11 \}} \alpha_q \log q < C(d) \cdot 1.05^2 \log^2 n \log y,
\end{equation}
where 
$$
C(d) = 2 \sum_{q \in \{ 3, 5, 7, 11 \}} c(d,q) \log q.
$$
We have
$$
C(7) = C(15)=C(55)=2 \left( \frac{576 \log 2}{\log^3 3} + \frac{480}{\log^2 5} +\frac{672}{\log^2 7} + \frac{1056}{\log^2 11} \right) < 1696
$$
and
$$
C(231) = 2 \left(  \frac{864 \log 2}{\log^3 3}  + \frac{1440 \log 2}{\log^3 5} +\frac{2016 \log 2}{\log^3 7} + \frac{1056}{\log^2 11} \right) < 2129.
$$

From (\ref{eqn:fracprincipal}) and (\ref{lambda!}),  we can write
\begin{equation} \label{vch}
\Lambda = n \log \left( \tau \delta \right) + b_2 \log \left( \gamma \right) + j \pi i,
\end{equation}
 with $b_2$ as in (\ref{bee2}), while, if $d=231$, we also have
\begin{equation} \label{vch3}
\Lambda^\prime = 3 \Lambda = n \log \left( \tau^\prime \delta^3 \gamma^{\epsilon_n} \right) + 2 \log \left( \gamma \right) + j^\prime \pi i.
\end{equation}
In each case, we take the principal branches of the logarithms, choose  $\tau, \tau^\prime\in \{ -1, 1 \}$ so that
$\mbox{Im} ( \log \left( \tau \delta \right) )$ and 
$ \mbox{Im}  (\log \left( \tau^\prime \delta^3 \gamma^{\epsilon_n} \right) )$ have opposite signs to $\mbox{Im}(\log \gamma)$, and take integers  $j$ and $j^\prime$  so that 
$|\Lambda|$ and $|\Lambda^\prime|$ are minimal. 
Notice that, with these choices, 
\begin{equation} \label{righteous1}
n  |\log (\tau \delta) | =
\vert \log\left( \gamma \right) \vert + |j| \,\pi  \pm  \vert\Lambda \vert, \; \mbox{ if } \; d \in \{ 7, 15, 55 \}
\end{equation}
and
\begin{equation} \label{righteous2}
n  |\log (\tau^\prime \delta^3 \gamma^{\epsilon_n}) | =
2 \,\vert \log\left( \gamma \right) \vert + |j^\prime| \,\pi  \pm  \vert\Lambda^\prime \vert.
\end{equation}
Note further that we have 
$$
\begin{array}{|c|c|} \hline
d & |\log \left(  \gamma\right)| \\ \hline
7  & \arccos(1/8) \\
15  & \arccos(7/8)  \\
55  & \arccos(3/8)  \\
231  & \arccos(5/16).  \\ \hline
 \end{array}
 $$

Assume first that inequality (\ref{ass}) fails to hold. Then, from (\ref{whopper}), we have
$$
n  < \frac{2 \log 10}{\log y}+ C(d) \cdot  1.05^2 \log^2 n,
$$
contradicting Lemma \ref{ybound}, (\ref{starter!}) and $C(d)<2129$. 
It follows that inequality (\ref{ass})  holds and hence we may conclude, from Lemma \ref{smallform},  that
$$
\log \left|  \Lambda \right| < 0.75 +  \frac{1}{2} C(d) \cdot  1.05^2 \log^2 n \log y - \frac{n}{2} \log y.
$$
From Lemma \ref{ybound}, (\ref{starter!}) and $C(d) < 2129$, we find, in all cases, that
\begin{equation} \label{upsies}
\log \left|  \Lambda \right| < -0.499  \,  n \,  \log y.
\end{equation}

\subsubsection{Linear forms in three logarithms}

To deduce an initial lower bound upon the linear form in logarithms $\left| \Lambda \right| $, we will use 
the following,  the main result (Theorem 2.1) of Matveev \cite{Mat}.
\begin{thm}[Matveev] \label{Matveev} 
Let $\mathbb{K}$ be an algebraic number field of degree $D$ over $\mathbb{Q}$ and put $\chi=1$ if $\mathbb{K}$ is real, $\chi=2$ otherwise. Suppose that $\alpha_1, \alpha_2, \ldots, \alpha_{n_0} \in \mathbb{K}^*$ with absolute logarithmic heights $h(\alpha_i)$ for $1 \leq i \leq n_0$, and suppose that
$$
A_i \geq \max \{ D \, h (\alpha_i), \left| \log \alpha_i \right| \}, \; 1 \leq i \leq n_0,
$$
for some fixed choice of the logarithm. Define
$$
\Lambda = b_1 \log \alpha_1 + \cdots + b_{n_0} \log \alpha_{n_0},
$$
where the $b_i$ are integers and  set
$$
B = \max \{ 1, \max \{ |b_i| A_i/A_{n_0} \; : \; 1 \leq i \leq n_0 \} \}.
$$
Define, with $e := \exp(1)$, further, 
$$
\Omega =A_1 \cdots A_{n_0}, 
$$
$$
C(n_0) = C(n_0,\chi) = \frac{16}{n_0! \chi} e^{n_0} (2n_0+1+2 \chi) (n_0+2)(4n_0+4)^{n_0+1} \left( en_0/2 \right)^{\chi},
$$
$$
C_0 = \log \left( e^{4.4 n_0+7} n_0^{5.5} D^2 \log ( e D) \right) \; \mbox{ and } \; W_0 = \log \left(
1.5 e B D \log (eD) \right).
$$
Then, if $\log \alpha_1, \ldots, \log \alpha_{n_0}$ are linearly independent over $\mathbb{Z}$ and $b_{n_0} \neq 0$, we have
$$
\log \left| \Lambda \right| > - C(n_0) \, C_0 \, W_0 \, D^2 \, \Omega.
$$
\end{thm}

We apply Theorem \ref{Matveev} to $\Lambda$ as given in (\ref{vch}), with 
$$
D=2, \; \chi = 2, \; n_0=3, \; b_3=n, \; \alpha_3 =\tau \delta, \;\alpha_2 = \gamma, \; b_1=j, \; \alpha_1=-1,
$$
and $b_2$ as in (\ref{bee2}). 

We may thus take  
$$
A_3 =  \log y,  \; A_2 = 3 \log 2, \; A_1 = \pi \;  \mbox{ and } \;  B =n.
$$
Since
$$
4 \, C(3) \, C_0 = 2^{18} \cdot 3 \cdot 5 \cdot 11 \cdot e^5  \cdot \log \left( e^{20.2}  \cdot 3^{5.5} \cdot 4 \log (2e) \right) < 
1.80741 \times 10^{11},
$$
and 
$$
W_0 =  \log \left(  3 e n \log (2e) \right) < 2.63+\log n,
$$
we may  therefore conclude that
$$
\log  \left| \Lambda \right| > - 1.181 \times 10^{12}  \left(  2.63+\log n \right) \log y.
$$
It thus follows from (\ref{upsies}) that
$$
n < 2.37 \times 10^{12} (\log n+2.63),
$$
whence
\begin{equation} \label{gumby}
n <8.22 \times 10^{13}.
\end{equation}

To improve this inequality, we appeal to a sharper but less convenient lower bound for linear forms in three complex logarithms, due to Mignotte  (Theorem 2 of \cite{Mig2}).

\begin{thm}[Mignotte] \label{miggy}
Consider three non-zero  algebraic numbers $\alpha_1$, $\alpha_2$
and $\alpha_3$, which are either all real and ${}>1,$ or all complex of modulus
one and all ${}\not=1$. Further, assume that the three numbers $\alpha_1, \alpha_2$ and $\alpha_3$ are either all multiplicatively independent, or that two of the numbers are multiplicatively independent and the third is a root of unity.
We also consider three positive
rational integers $b_1$, $b_2$, $b_3$ with $\gcd(b_1,b_2,b_3)=1$, and the linear form
$$
   \Lambda = b_2\log \alpha_2-b_1 \log \alpha_1-b_3\log \alpha_3 ,
$$
where the logarithms of the $\alpha_i$ are arbitrary determinations of the logarithm,
but which are all real or all purely imaginary.
We assume that
$$
0 < |\Lambda| < 2\pi/w,
$$
where $w$ is the maximal order of a root of unity in $\mathbb{Q}(\alpha_1,\alpha_2,\alpha_3)$. Suppose further that
\begin{equation} \label{frisky}
   b_2 |\log \alpha_2| =
 b_1\,\vert \log \alpha_1 \vert+  b_3 \,\vert\log \alpha_3\vert \pm  \vert\Lambda\vert 
\end{equation}
and put
$$
d_1 = \gcd(b_1,b_2), \; \;   d_3 = \gcd(b_3,b_2) \; \mbox{ and } \; b_2=d_1 b_2^\prime = d_3 b_2^{\prime\prime} 
$$
Let
$K, L, R, R_1, R_2, R_3, S, S_1, S_2, S_3, T, T_1, T_2, T_3$
be positive rational integers with
$$
K \geq 3, \; L \geq 5, \; R > R_1+R_2+R_3, \; S > S_1+S_2+S_3 \; \mbox{ and } \; T > T_1+T_2+T_3
$$

Let $\rho\ge 2$  be a real number. Let $a_1, a_2$ and $a_3$ be real numbers such that
$$
   a_i \ge  \rho \vert \log \alpha_i \vert
   - \log  \vert \alpha_i\vert +2 D \,{\rm h}\kern .5pt(\alpha_i), \qquad
   i \in \{1, 2, 3 \},
$$
where
$\,D=[\mathbb{Q}(\alpha_1,\alpha_2,\alpha_3) : \mathbb{Q}]\bigm/[\mathbb{R}(\alpha_1,\alpha_2,\alpha_3) : \mathbb{R}]$, and set
$$
U = \left( \frac{KL}{2}+\frac{L}{4}-1- \frac{2K}{3L} \right) \log \rho.
$$
Assume further that
\begin{equation} \label{needed}
U \geq (D+1) \log (K^2L) + gL(a_1R+a_2S+a_3T) + D(K-1) \log b - 2 \log (e/2),
\end{equation}
where 
$$
g=\frac{1}{4}-\frac{K^2L}{12RST} \; \mbox{ and } \;  b=(b_2^\prime \eta_0)(b_2^{\prime\prime} \zeta_0) \left( \prod_{k=1}^{K-1} k! \right)^{-\frac{4}{K(K-1)}},
$$ 
with
$$
\eta_0=\frac{R-1}{2}+\frac{(S-1)b_1}{2b_2} \; \mbox{ and } \; \zeta_0=\frac{T-1}{2} + \frac{(S-1)b_3}{2b_2}.
$$
Put 
$$
\mathcal{V} = \sqrt{(R_1+1)(S_1+1)(T_1+1)}.
$$
If, for some positive real number $\chi$, we have
\begin{enumerate}[label=(\roman*)]
\item $(R_1+1)(S_1+1)(T_1+1) > K \mathcal{M}$, \\
\item $\mbox{Card} \{ \alpha_1^r \alpha_2^s\alpha_3^t \; : \; 0 \leq r \leq R_1, \; 0 \leq s \leq S_1, \; 0 \leq t \leq T_1 \} > L$, \\
\item $(R_2+1)(S_2+1)(T_2+1) > 2 K^2$, \\
\item $\mbox{Card} \{ \alpha_1^r \alpha_2^s\alpha_3^t \; : \; 0 \leq r \leq R_2, \; 0 \leq s \leq S_2, \; 0 \leq t \leq T_2 \} > 2K L$,  and \\
\item $(R_3+1)(S_3+1)(T_3+1) > 6 K^2 L$, \\
\end{enumerate}
where
$$
\mathcal{M} =  \max\Bigl\{R_1+S_1+1,\,S_1+T_1+1,\,R_1+T_1+1,\,\chi \; \mathcal{V} \Bigr\}, \; \; 
$$
then either
\begin{equation} \label{smosh}
\left| \Lambda \right| \cdot \frac{LSe^{LS|\Lambda|/(2b_2)}}{2|b_2|} > \rho^{-KL},
\end{equation}
or at least one of the following conditions ({\bf{C1}}), ({\bf{C2}}), ({\bf{C3}}) holds :

\vskip2ex
\noindent ({\bf{C1}})  \hskip4ex $|b_1| \leq R_1$  and $|b_2| \leq S_1$ and  $|b_3| \leq T_1$,

\vskip1.5ex
\noindent ({\bf{C2}})  \hskip4ex $|b_1| \leq R_2$  and $|b_2| \leq S_2$ and  $|b_3| \leq T_2$,

\vskip1.5ex
\noindent \mbox{({\bf{C3}})}   {\bf either} there exist non-zero rational integers $r_0$ and $s_0$ such that
\begin{equation} \label{rups2}
   r_0b_2=s_0b_1
\end{equation}
with
\begin{equation} \label{rups3}
   |r_0|
   \le \frac{(R_1+1)(T_1+1)}{\mathcal{M}-T_1}
    \;  \mbox{ and } \; 
   |s_0| 
   \le \frac{(S_1+1)(T_1+1)}{\mathcal{M}-T_1},
\end{equation}
{\bf or}
there exist rational integers  $r_1$, $s_1$, $t_1$ and $t_2$, with
$r_1s_1\not=0$, such that
\begin{equation} \label{rups4}
   (t_1b_1+r_1b_3)s_1=r_1b_2t_2, \qquad \gcd(r_1, t_1)=\gcd(s_1,t_2 )=1,
\end{equation}
which also satisfy
$$
    |r_1s_1|
   \le \gcd(r_1,s_1) \cdot
   \frac{(R_1+1)(S_1+1)}{\mathcal{M}-\max \{ R_1, S_1 \}},
$$
$$
   |s_1t_1| \le \gcd(r_1,s_1) \cdot
   \frac{(S_1+1)(T_1+1)}{\mathcal{M}-\max \{ S_1, T_1 \}} 
  $$
and
$$
   |r_1t_2| % &
   \le \gcd(r_1,s_1) \cdot
  \frac{(R_1+1)(T_1+1)}{\mathcal{M}- \max \{ R_1, T_1 \}}.
$$
Moreover, when $t_1=0$ we can take $r_1=1$, and
when $t_2=0$ we can take $s_1=1$.
\end{thm}

We will apply this result to our $\Lambda$ (if $d \in \{ 7, 15, 55 \}$) or $\Lambda^\prime$ (if $d=231$). 
To do this, we must distinguish between a number of cases, depending on $d$ and the signs of the coefficients in (\ref{vch}) or  (\ref{vch3}).
By way of example, suppose first that $d=231$. If we have $j^\prime=0$ or $j^\prime = \pm n$, then $\Lambda^\prime$ reduces to a linear form in two logarithms and we may appeal to Corollary 1 of Laurent \cite{Lau}; actually, what we state here is specialized for our purposes and follows from the arguments of \cite{Lau} (see pages 346 and 347) after a short computation (in each case using, in the notation of \cite{Lau}, values of $\mu$ in $[0.555,0.562]$ and $\rho$ in $[6.12,6.31]$).
\begin{thm}[Laurent] \label{LFL2}
Consider the linear form 
$$
\Lambda = c_2 \log  \beta_2 - c_1 \log \beta_1, 
$$
where $c_1$ and $c_2$ are positive integers, and $\beta_1$ and $\beta_2$ are multiplicatively independent algebraic numbers. Define
$\,D=[\mathbb{Q}(\beta_1,\beta_2) : \mathbb{Q}]\bigm/[\mathbb{R}(\beta_1,\beta_2) : \mathbb{R}]$ and set
$$
b^\prime =\frac{c_1}{D \log B_2}+ \frac{c_2}{D \log B_1},
$$
where $B_1, B_2 > 1$ are real numbers such that
$$
\log B_i \geq \max \{ h (\beta_i), |\log \beta_i|/D, 1/D \}, \; \; i \in \{ 1, 2 \}.
$$
Then
$$
\log \left| \Lambda \right| \geq -C D^4 \left( \max \{ \log b^\prime + 0.21, m_1/D, 1 \} \right)^2 \log B_1 \log B_2,
$$
for each pair $(m_1,C)$ in the following set
$$
\begin{array}{c}
 \left\{ (14,28.161), (14.5,27.812), (15,27.486), (15.5,27.182), (16,26.896), (16.5,26.627), (17,26.374),    \right. \\
\left. (17.5,26.136), (18,25.911),  (18.5,25.697), (19,25.495), (19.5,25.303), (20,25.120)  \right\}.\\
\end{array}
$$
\end{thm}

If $j^\prime = 0$, we apply this with
$$
c_2=n, \; \beta_2 =  \tau^\prime \delta^3 \gamma^{\epsilon_n}, \; c_1=2, \; \beta_1 = 1/\gamma, \; D=1, 
$$
whence
$$
h(\beta_2) \leq \frac{3}{2} \log (y), \;  \; h(\beta_1)= \frac{3 \log 2}{2}, \;
$$
and we can take
$$
\log B_2=\frac{3}{2} \log (y) \; \mbox{ and } \; \log B_1 = \arccos( 5/16).
$$
Choosing $(m_1,C)=(18,25.911)$, it follows that
$$
\log \left| \Lambda^\prime \right| \geq -39 \arccos( 5/16) \left( \log (n)+ 0.18 \right)^2 \log y,
$$
whereby, from (\ref{upsies}) and $\Lambda^\prime = 3 \Lambda$, $n < 8300$, contradicting (\ref{starter!}). Similarly, if $j^\prime = \pm n$, we may apply Theorem \ref{LFL2} with 
$$
c_2=n, \; c_1=2, \; \beta_1 = 1/\gamma, \; D=1, 
$$
where
$$
\log \beta_2 = \log ( \tau^\prime \delta^3 \gamma^{\epsilon_n} ) \pm \pi i.
$$
We can once again choose
$$
\log B_2=\frac{3}{2} \log (y) \; \mbox{ and } \; \log B_1 = \arccos( 5/16),
$$
and derive a contradiction from taking $(m_1,C)=(18,25.911)$. Note that from  (\ref{vch3}),  
$$
|j^\prime| \pi <  \pi  n + 2 \arccos(5/16) + 3 \cdot y^{-0.499 n} < \pi  n+ 2.51,
$$
whereby $|j^\prime| \leq n$. We may thus suppose that $|j^\prime| < n$ and $j^\prime \neq 0$ (so that, in particular, we have $\gcd (j^\prime, n)=1$).

We will now apply Theorem \ref{miggy}.
From (\ref{righteous2}), 
we can take, in the notation of Theorem \ref{miggy} and  writing $\upsilon = -j^\prime/|j^\prime|$,
\begin{equation} \label{case1}
b_1=2, \; \alpha_1 = \gamma^{-\upsilon} , \; b_2=n, \; \alpha_2 =(\tau^\prime \delta^3 \gamma^{\epsilon_n})^\upsilon, \;  b_3=|j^\prime|\; \mbox{ and }  \;  \alpha_3=-1.
\end{equation}
It follows that
$$
h(\alpha_1)= \frac{3 \log 2}{2}, \; h(\alpha_2) \leq \frac{3}{2} \log (y) \; \mbox{ and } \;  h(\alpha_3) =0.
$$
We can thus choose 
$$
a_1=\rho \arccos(5/16) + 3\log (2), \; a_2=  \rho  \pi + 3 \log (y) \; \mbox{ and } \; a_3 = \rho \pi.
$$
As noted in \cite{BMS2}, if we suppose that $m \geq 1$ and define
\begin{equation} \label{excellent1}
\begin{array}{c}
K=[m L a_1 a_2 a_3], 
\; R_1 = [c_1a_2a_3], \; S_1=[c_1a_1a_3], \; T_1=[c_1a_1a_2],  \; R_2 = [c_2a_2a_3], \\
\\
S_2=[c_2a_1a_3], \; T_2=[c_2a_1a_2], \;
R_3 = [c_3a_2a_3], \; S_3=[c_3a_1a_3] \; \mbox{ and } \; T_3=[c_3a_1a_2],  \\
\end{array}
\end{equation}
where
\begin{equation} \label{excellent2}
\begin{array}{c}
c_1=\max \{ (\chi m L)^{2/3}, (2mL/a_1)^{1/2} \}, \; c_2=\max \{ 2^{1/3} (m L)^{2/3}, (m/a_1)^{1/2} L \}  \\ \\
\mbox{ and } c_3=(6m^2)^{1/3} L, \\
\end{array}
\end{equation}
then conditions (i)-(v) are automatically satisfied. It remains to verify inequality (\ref{needed}).

Define
$$
R = R_1+R_2+R_3+1, \; S = S_1+S_2+S_3+1 \; \mbox{ and } \; T = T_1+T_2+T_3+1.
$$
We choose 
$$
\rho = 5.9, \; L=206, \; m=25 \; \mbox{ and } \; \chi=2.89,
$$
so that
$$
c_1 =  (\chi m L)^{2/3}, \; \; c_2=2^{1/3} (m L)^{2/3},
$$
and we have
$$
K = [K_1+K_2\log (y)],
$$
where
$$
K_1=16759141.618\ldots   \; \mbox{ and } \; K_2=2712508.708\ldots.
$$
We thus have
$$
S_1=106229, \; S_2=65966 \; \mbox{ and } \; S_3= 561893.   
$$
Since Lemma \ref{ybound}, (\ref{starter!}) and $N(231) =1.2 \times 10^9$ together imply that
\begin{equation} \label{wanda}
\log y > 22.2,
\end{equation} 
we find, after a little work, that $\mathcal{M} = \chi  \mathcal{V}$ and that $g < 0.2438$.  

Since $\gcd (j^\prime, n)=1$, we have
$$
d_1=d_3=1, \; b_2^\prime = b_2^{\prime\prime} =  n, \; 
$$
and it follows  that
$$
\eta_0 = \frac{1}{2} \left( R_1+R_2+R_3 \right) + \frac{1}{n} \left( S_1+S_2+S_3 \right) <  718258+116252 \log y,
$$
and, from $|j^\prime| < n$, 
$$
\zeta_0 = \frac{1}{2} \left( T_1+T_2+T_3 \right) + \frac{-j^\prime}{2n} \left( S_1+S_2+S_3 \right)  < \frac{1}{2} \left( T_1+T_2+T_3 + S_1+S_2+S_3 \right),
$$
whereby
$$
\zeta_0  < 734089 + 59408 \log y.    
$$
From Lemma 3.4 of \cite{Mig2}, we have the inequality
\begin{equation} \label{mess}
\log \left( \prod_{k=1}^{K-1} k! \right)^{\frac{4}{K(K-1)}} \geq 2 \log K -3+\frac{2 \log \left( 2 \pi K/e^{3/2} \right)}{K-1}-
\frac{2+6 \pi^{-2}+\log K}{3K(K-1)},
\end{equation}
whence, from $K > 10^6$,
$$
\log \left( \prod_{k=1}^{K-1} k! \right)^{\frac{4}{K(K-1)}} > 2 \log K -3.
$$
It follows, from (\ref{wanda}), that 
$$
b < e^3 n^2 \frac{\left( 718258+116252 \log y \right) \left( 734089 + 59408 \log y \right) }{\left( 16759141.6+2712508.7 \log y \right)^2}
< 0.023062 n^2 < 1.559 \times 10^{26},
$$
where the last inequality is a consequence of (\ref{gumby}). 
The right-hand-side of inequality (\ref{needed}) is thus bounded above by 
$$
4 \log (K)  +2.051 \times 10^9 + 3.319 \times 10^8  \log (y) + 60.312 K 
$$
while the left-hand-side satisfies
$$
U >  182.814 K + 89.635.
$$
If inequality (\ref{needed}) fails to hold, it follows that
$$
122.502 K < 4 \log (K)  + 2.051 \times 10^9 + 3.318 \times 10^8  \log (y),
$$
contradicting 
$$
K>  16759141 + 2712508 \log y 
$$
and (\ref{wanda}).

Note that we have
$$
\frac{LSe^{LS|\Lambda|/(2b_2)}}{2|b_2|} = \frac{75611167 \, e^{75611167 |\Lambda|/n}}{n} 
%<  \frac{ 22292930  \mbox{exp} \left(\frac{22292930 }{n y^{0.4993}}}}{n} \right)}
$$
and hence, from (\ref{upsies}),
$$
\frac{LSe^{LS|\Lambda|/(2b_2)}}{2|b_2|} 
<  \frac{ 75611167 \,  \mbox{exp} \left( \frac{75611167}{n y^{0.499n}}\right)}{n} < 0.127,
$$
where the last inequality is a consequence of Lemma \ref{ybound} and (\ref{starter!}).
If we have inequality (\ref{smosh}), it thus follows that
$$
\log \left| \Lambda \right| >  2 -365.65 K.
$$
Once again appealing to (\ref{upsies}), we find that
$$
0.499 \, n \log y < 365.65 K- 2 < 365.65 \left( 16759141.62+2712508.71 \log y \right)
$$
and so
$$
n < 1.988 \times 10^9 + \frac{1.229 \times 10^{10}}{\log y},
$$
whence, from  (\ref{wanda}),
\begin{equation} \label{starter2}
n < 2.6 \times 10^9.
\end{equation}

If, on the other hand, inequality (\ref{smosh}) fails to be satisfied, from inequality (\ref{starter!}) and our choices of $S_1$ and $S_2$, necessarily ({\bf{C3}}) holds. If (\ref{rups2}) holds then $n \mid s_0$, where 
$$
|s_0| \leq \frac{(S_1+1)(T_1+1)}{\mathcal{M}-T_1} \leq  \frac{(S_1+1)(T_1+1)}{R_1+1} < \frac{106230 \, (106230.73 + 17193.55 \log y )}{207878 + 33645 \log y} < 54287, 
$$
from calculus.
It follows that necessarily $s_0=0$,  a contradiction. We thus have (\ref{rups4}). In particular,
\begin{equation} \label{two-log}
  \left( 2 t_1 +r_1 |j^\prime| \right) s_1=r_1t_2 n,
\end{equation}
for integers $r_1, s_1, t_1, t_2$ with $\gcd(r_1, t_1)=\gcd(s_1,t_2 )=1$,
\begin{equation} \label{ooze}
   |s_1t_1| \le  \gcd(r_1,s_1) \cdot
   \frac{(S_1+1)(T_1+1)}{\chi \mathcal{V}-T_1} <  \gcd(r_1,s_1) \cdot 81
\end{equation}
and
\begin{equation} \label{ooze2}
   |r_1s_1| \le \gcd(r_1,s_1) \cdot
   \frac{(R_1+1)(S_1+1)}{\chi \mathcal{V} -  R_1 }
   <   \gcd(r_1,s_1) \cdot 158,
\end{equation}
again via calculus.
It follows that
\begin{equation} \label{tbound}
|t_1| \leq 80 \; \; \mbox{ and } \; \; |r_1| \leq 157.
\end{equation}
Since $r_1$ is coprime to $t_1$, necessarily $r_1 \mid 2s_1$, while the fact that $\gcd (s_1,t_2)=1$ while  $n > 1.2 \times 10^9$ is prime, together imply that $s_1 \mid r_1$. 
We thus have that $r_1 = \pm s_1$ or $r_1=\pm 2 s_1$, whence, from (\ref{two-log}),  
\begin{equation} \label{two-logA}
r_1 j^\prime = t_3 n \pm 2 t_1,
\end{equation}
where $t_3 = \pm t_2$ or $t_3=\pm t_2/2$.
We can thus rewrite $r_1\Lambda^\prime$ as a linear form in two logarithms,
$$
r_1 \Lambda^\prime = n \log \alpha - 2 \log \beta,
$$
where
$$
\log \alpha = r_1  \log \left( \tau^\prime \delta^3 \gamma^{\epsilon_n} \right) +t_3 \pi i \; \mbox{ and } \; 
\log \beta = -r_1 \log (\gamma) \pm  t_1  \pi i.
$$
We apply Theorem \ref{LFL2} with 
$$
D=1, \; c_2=n, \; \beta_2 =  \alpha, \; c_1 =2 \; \mbox{ and } \; \log \beta_1 = \log \beta.
$$
We may take
$$
\log B_2=\frac{3}{2} |r_1| \log (y) \; \mbox{ and } \; \log B_1 =
|r_1| \arccos( 5/16) +|t_1| \pi,
$$
whence
$$
b^\prime= \frac{n}{|r_1| \arccos( 5/16) +|t_1| \pi} + \frac{4}{3 |r_1| \log (y) }.
$$
From (\ref{tbound}), we thus have
$$
\log B_1  < 448.05.
$$
Again choosing $(m_1,C)=(18,25.911)$, we conclude that
$$
\log | \Lambda^\prime| >  -39 \times 157 \times 448.05 \log ^2 (n)  \log (y) - \log 157.
$$
From (\ref{upsies}) and $\Lambda^\prime = 3 \Lambda$, we thus have
$$
0.499  \,  n  <  \frac{\log 471}{\log y} + 2.744 \times 10^6  \log ^2 (n),
$$
whence, from (\ref{wanda}), we once again obtain inequality
(\ref{starter2}).

We now iterate this argument, using slightly more care and, in each case, assuming that  $n > N(231)=1.2 \times 10^9$. We begin by taking
$$
\rho = 5.8, \; L= 144, \; m=27 \; \mbox{ and } \; \; \chi = 3.14, 
$$
and, arguing as previously, find that either
\begin{equation} \label{starter3}
n < 1.3 \times 10^9,
\end{equation}
or that we have (\ref{two-logA}) with 
\begin{equation} \label{tbound2}
|t_1| \leq 68 \; \; \mbox{ and } \; \; |r_1| \leq 133.
\end{equation}
From (\ref{starter2}), we have that
$$
b^\prime= \frac{n}{|r_1| \arccos( 5/16) +|t_1| \pi} + \frac{4}{3 |r_1| \log (y) } < e^{18-0.21},
$$
provided, crudely, $|r_1| \geq 39$. Applying  Theorem \ref{LFL2} as previously, once again with $(m_1,C)=(18,25.911)$,
if $39 \leq |r_1| \leq 133$, we find that
$$
\log | \Lambda^\prime| >  -39 \times 133 \times 380.28 \times 18^2  \log (y) - \log 133,
$$
while, in case $1 \leq |r_1| \leq 38$, 
\begin {equation} \label{fool}
\log | \Lambda^\prime| >  -39 \times 38 \times 261.25 \log ^2 (n)  \log (y) - \log 38.
\end{equation}
Once more, we have, in either case, inequality \ref{starter3}.

To finish the case $d=231$, we repeat this argument, only now using
$$
\rho = 6.3, \; L= 125, \; m=27 \; \mbox{ and } \; \; \chi = 3.28, 
$$
and appealing to Theorem \ref{LFL2}  with $(m_1,C)=(17,26.38)$. We conclude then, in all cases with $d=231$, that
$$
n < 1.2 \times 10^9 = N(231).
$$

We argue similarly for $d \in \{ 7, 15, 55 \}$, applying Theorem \ref{miggy} to $\Lambda$ as in (\ref{vch}). As in the case $d=231$, we either have inequality (\ref{smosh}), or we find ourselves in the degenerate case where, analogous to (\ref{two-logA}), we deduce the existence of ``small'' integers $t_1, t_3$ and $r_1$, the last nonzero, such that $r_1 j = t_3 n \pm 2 t_1$.
In the latter case, we  rewrite  $r_1\Lambda$ as a linear form in two logarithms,
$$
r_1 \Lambda= n \log \alpha -  \log \beta,
$$
where
$$
\log \alpha = r_1  \log \left( \tau \delta \right) +t_3 \pi i \; \mbox{ and } \; 
\log \beta = -r_1 \log (\gamma) \pm  2t_1  \pi i,
$$
and apply Theorem \ref{LFL2} with 
$$
D=1, \; c_2=n, \; \beta_2 =  \alpha, \; c_1 =1, \; \log \beta_1 = \log \beta, \;  \log B_2=\frac{1}{2} |r_1| \log (y), \;  \log B_1 =
|r_1| |\log \gamma| +2 |t_1| \pi
$$
and
$$
b^\prime= \frac{n}{|r_1| |\log\gamma|  +2 |t_1| \pi} + \frac{2}{|r_1| \log (y) }.
$$
We make our parameter choices for Theorems \ref{miggy} and \ref{LFL2} as in the following table. In practice, once we have a reasonable upper bound upon $n, |r_1|$ and $|t_1|$, we loop over $|r_1|$ and $|t_1|$, compute upper and lower bounds upon $b^\prime$ and, in case $b^\prime$ potentially exceeds $\exp(m_1-0.21)$, deduce sharpened versions of inequality (\ref{fool}).

$$
\begin{array}{|c|c|c|c|c|c|c|} \hline
d & \rho & L & m & \chi & (m_1,C) & \mbox{upper bound upon $n$} \\ \hline
7 & 4.8 & 205 & 45 & 2.96 & (20,25.120) & 1.27 \times 10^9 \\
7 & 5.1 & 142 & 40 & 3.91 & (20,25.120) & 6.44 \times 10^8 \\
7 & 6.2 & 100 & 40 & 3.00 & (14,28.161) & 6 \times 10^8 = N(7) \\ \hline
15 & 5.8 & 180 & 36 & 2.77 &  (20,25.120) & 6.00 \times 10^8 \\
15 & 5.1 & 140 & 40 & 3.00 & (20,25.120)  & 4 \times 10^8 = N(15) \\ \hline
55 & 5.7 & 189 & 33 & 3.00 & (20,25.120)  & 1.11 \times 10^9 \\
55 & 5.2 & 144 & 35 & 3.41 &  (18,25.911) & 5.22 \times 10^8 \\
55 & 6.2 & 100 & 40 & 3.00 & (14,28.161) & 5 \times 10^8 = N(55) \\
 \hline
\end{array}
$$

\vskip1.5ex
\noindent This, with Proposition \ref{prop:yeven}, completes the proof of Theorem \ref{Main-thm}.

As mentioned previously, of note here is that the bounds we obtain upon the exponent $n$ for the equation
$x^2+c^2d=y^n$,
with $d \in \{ 7, 15, 55, 231 \}$ and $c$ an $S$-unit, $S=\{ 3, 5, 7, 11 \}$,
are essentially identical to those deduced for the  simpler equation
$x^2+d=y^n$.
This is admittedly not immediately apparent from perusal of Section 15 of \cite{BMS2}, where the treatment of the degenerate cases (which reduce to linear forms in two logarithms) requires some modification.

%---------------------------------------------------------
\section{Concluding remarks}
%---------------------------------------------------------

Extending the results of this paper to the more general equation 
$$
x^2+D = y^n, \; \; \gcd (x,y)=1, \; \; D > 0, \; \; P(D) \leq 13
$$
is probably computationally feasible with current technology, if one is suitably enthusiastic, while the equation 
$$
x^2+D = y^n, \; \; \gcd (x,y)=1, \; \; D > 0, \; \; P(D) \leq 17
$$
is certainly out of reach without the introduction of fundamentally new ideas. The two main obstructions arise from both large exponents $n$ (where the corresponding spaces of modular forms have extremely large dimensions) and moderately small ones (where one will encounter Thue-Mahler equations with very many associated $S$-unit equations). 
Additionally, it is possible to relax the restriction that $\gcd (x,y)=1$ in (\ref{eq-main}) (at least provided this $\gcd$ is odd), though the computational difficulties increase substantially since, once again, the spaces of modular forms one encounters have significantly higher dimensions.

%---------------------------------------------------------
\section{Acknowledgements}
%---------------------------------------------------------

The authors would like to thank the anonymous referee, Vandita Patel and Pedro-Jose Carzola Garcia for numerous helpful suggestions that improved both the exposition and correctness of this paper.

%---------------------------------------------------------

\bibliographystyle{abbrv}
\bibliography{samir}
\end{document}